\documentclass{amsart}
\usepackage{amsfonts}
\usepackage{amsmath,amssymb}
\usepackage{amsthm}
\usepackage{amscd}
\usepackage{graphics}
\usepackage{graphicx}
{
\theoremstyle{definition}
\newtheorem{Def}{Definition}
\newtheorem{Ex}{Example}
\newtheorem{Rem}{Remark}

}
\newtheorem{Cor}{Corollary}
\newtheorem{Prop}{Proposition}
\newtheorem{Thm}{Theorem}

\begin{document}
\title{Constructions of round fold maps on $C^{\infty}$ bundles}
\author{Naoki Kitazawa}
\address{Department of Mathematics, Tokyo Institute of Technology, 
2-12-1 Ookayama, Meguro-ku, Tokyo 152-8551, JAPAN, Tel +81-(0)3-5734-2205, 
Fax +81-(0)3-5734-2738,
}
\email{kitazawa.n.aa@m.titech.ac.jp}
\maketitle
\begin{abstract}
In this paper, we construct {\it round fold maps} or {\it stable fold maps} with concentric singular value sets introduced
 by the author \cite{kitazawa 2} on $C^{\infty}$ bundles over spheres or bundles over more general manifolds. The class of round fold maps
 includes {\it special generic maps} on spheres (see \cite{saeki 2} for example) and such maps have been 
 constructed on $C^{\infty}$ bundles over standard spheres of dimensions larger than $1$ and connected
 sums of $C^{\infty}$ bundles over standard spheres of dimensions larger than $1$ whose fibers are standard spheres, for example, in previous studies by the
 author (\cite{kitazawa 4}, \cite{kitazawa 5}). In this paper, we obtain round fold maps
 and the diffeomorphism types of their source manifolds which do not appear in these studies in new manners.
\end{abstract}
\section{Introduction}
\label{sec:1}
{\it Fold maps} are important in generalizing the theory of Morse functions. Studies of
 such maps were started by Whitney (\cite{whitney}) and Thom (\cite{thom}) in the 1950's. A {\it fold map} is a $C^{\infty}$ map
 whose singular points are of the form
$$(x_1, \cdots, x_m) \mapsto (x_1,\cdots,x_{n-1},\sum_{k=n}^{m-i}{x_k}^2-\sum_{k=m-i+1}^{m}{x_k}^2)$$ for two positive integers $m \geq n$ and an integer $0 \leq i \leq m-n+1$. A Morse function
 is a fold map of course.
  For a fold map from a closed $C^{\infty}$ manifold of dimension $m$ into a $C^{\infty}$ manifold of dimension $n$
 (without boundary), the followings hold where $m \geq n \geq 1$.
\begin{enumerate}
\item The {\it singular set}, or the set of all the singular points, is a closed $C^{\infty}$ submanifold of dimension {\rm (}$n-1${\rm )} of the source manifold. 
\item The restriction map to the singular set is a $C^{\infty}$ immersion of codimension $1$.
\end{enumerate}
We also note that if the restriction map to the singular set is a immersion with normal
 crossings, then it is {\it stable} (stable maps are important in the theory of global singularity; see \cite{golubitsky guillemin} for example). In \cite{kitazawa 2}, such
 a fold map is defined as a {\it stable} fold map. \\ 
\ \ \ Since around the 1990's, fold maps with additional conditions have been
 actively studied. For example, in \cite{burlet derham}, \cite{furuya porto}, \cite{saeki 2}, \cite{saeki sakuma} and \cite{sakuma}, {\it special
 generic} maps, which are fold maps whose singular points are of the form
$$(x_1, \cdots, x_m) \mapsto (x_1,\cdots,x_{n-1},\sum_{k=n}^{m}{x_k}^2)$$ for two positive integers $m \geq n$, were studied. In \cite{sakuma}, Sakuma studied
 {\it simple} fold maps, which are fold maps such that fibers of singular values do not have any connected component with more
 than one singular points (see also \cite{saeki}). For example, special generic maps are simple. In
 \cite{kobayashi saeki}, Kobayashi and Saeki
 investigated topology of {\it stable} maps including fold maps which are stable into the plane. In \cite{saeki suzuoka}, Saeki and Suzuoka
 found good properties of manifolds admitting stable maps whose regular fibers, or fibers of regular values, are disjoint unions of spheres. \\
\ \ \ Later, in \cite{kitazawa 2}, {\it round} fold maps, which will be mainly studied in this paper, were introduced. A {\it round} fold map is a fold
 map satisfying the followings.
\begin{enumerate}
\item The singular set is a disjoint union of standard spheres.
\item The restriction map to the singular set is a $C^{\infty}$ embedding.
\item The singular value set is a disjoint union of spheres embedded concentrically. 
\end{enumerate}
 Fold maps satisfying the definition of a round fold map have appeared in previous studies as the followings. 
\begin{itemize}
\item A lot of special generic maps on homotopy spheres. (See also \cite{saeki 2})
\item Fold maps represented as FIGURE 7 of \cite{kobayashi saeki} and Figure 8 of \cite{saeki suzuoka}.
\end{itemize}
 (Algebraic) topological properties of some round fold maps were studied in Theorem 1, Theorem 2 (homology groups of manifolds) and
 Theorem 3 (homotopy groups of manifolds) in \cite{kitazawa 2}. \\
\ \ \ Now, how about the homeomorphism and diffeomorphism types of the source manifolds? As answers, some examples of round fold maps
 and the diffeomorphism types of their source manifolds are given by the author in \cite{kitazawa}, \cite{kitazawa 4} and \cite{kitazawa 5}. For example, we have
 obtained round fold maps on $C^{\infty}$ bundles over standard spheres of dimensions larger than $1$ and manifolds represented as connected
 sums of $C^{\infty}$ bundles over standard spheres of dimensions larger than $1$ with fibers $C^{\infty}$ diffemorphic to standard
 spheres and spheres whose dimensions are positive. In this paper, as new answers to the question, we construct round fold maps of
 new types on closed $C^{\infty}$ manifolds having the
 structures of $C^{\infty}$ bundles over (exotic) spheres or bundles over more general manifolds. \\
\ \ \ This paper is organized as the following. \\ 
\ \ \ Section \ref{sec:2} is for preliminaries. We recall {\it fold maps} and introduce {\it stable} fold maps. We
 also recall {\it special generic} maps and {\it simple} fold maps. Finally we review the {\it Reeb space} of a smooth
 map, which is the space consisting of all the connected components of
 all the fibers of the smooth map. \\
\ \ \ In section \ref{sec:3}, we recall {\it round} fold maps and some terminologies
 on round fold maps such as {\it axes}, {\it proper cores}. We also recall a {\it topologically trivial} ({\it $C^{\infty}$ trivial}) round fold map. We introduce results
 on the diffeomorphism types of manifolds admitting topologically or $C^{\infty}$ trivial round fold maps shown by the author in \cite{kitazawa}, \cite{kitazawa 4} and \cite{kitazawa 5} (Proposition \ref{prop:4} and Proposition \ref{prop:5}). \\
\ \ \ In section \ref{sec:4}, we give some new examples of round fold maps with the diffeomorphism types of their source manifolds. Before constructions, we introduce a {\it locally $C^{\infty}$ trivial} round fold map as a map satisfying
 a kind of triviality around the connected components of the singular value set. Here, a $C^{\infty}$ trivial round fold map satisfies the definition of
 a locally $C^{\infty}$ trivial round fold map. Then, for example, we construct round fold maps on closed $C^{\infty}$ manifolds having the
 structures of trivial $C^{\infty}$ bundles over manifolds admitting locally $C^{\infty}$ trivial round fold maps (Proposition \ref{prop:7}) and closed $C^{\infty}$ manifolds having the
 structures of $C^{\infty}$ bundles over (exotic) spheres (Theorem \ref{thm:1} and Theorem \ref{thm:2}). \\
\ \ \ In section \ref{sec:5}, we define {\it P-operations}. {\it P-operations} are defined as generalizations of operations used to construct new round fold maps
 on manifolds having the structures of $C^{\infty}$ bundles over spheres in the previous section. We apply these operations
 to construct families of new round fold maps on $C^{\infty}$ $S^1$-bundles over manifolds admitting locally trivial $C^{\infty}$ round
 fold maps (Theorem \ref{thm:3}, Theorem \ref{thm:4} and Theorem \ref{thm:5}). \\
\ \ \ Throughout this paper, we assume that $M$ is a closed $C^{\infty}$ manifold of dimension $m$, that
 $N$ is a $C^{\infty}$ manifold of dimension $n$ without boundary, that $f:M \rightarrow N$ is a $C^{\infty}$ map and that $m \geq n \geq 1$.
 We denote the {\it singular set} of $f$, or the set consisting of all the singular points of $f$, by $S(f)$. \\

\section{Preliminaries}
\label{sec:2}
\subsection{Fold maps}
 First, we recall {\it fold maps}, which are simplest generalizations of Morse functions. For precise information, See \cite{golubitsky guillemin}, \cite{milnor} and \cite{milnor 2} and
 see also \cite{kitazawa 2}, \cite{kitazawa 3} and \cite{kitazawa 4}, for example.  

\begin{Def}
\label{def:1}
For a $C^{\infty}$ map $f:M \rightarrow N$, $p \in M$ is said to be a {\it fold} point of $f$ if at $p$ $f$ has
 the normal form $$f(x_1, \cdots, x_m):=(x_1,\cdots,x_{n-1},\sum_{k=n}^{m-i}{x_k}^2-\sum_{k=m-i+1}^{m}{x_k}^2)$$ and
 if all the singular points of $f$ are fold, then we call $f$ a {\it fold map}.
\end{Def}

If $p \in M$ is a fold point of $f$, then we can define $j:=\min \{i,m-n+1-i\}$ uniquely in the previous definition.
 We call $p$ a {\it fold point of index $j$} of $f$. We call a fold point of index $0$ a {\it definite} fold point of $f$ and
 we call $f$ a {\it special generic} map if all the singular points are definite fold. For special generic maps, see
 \cite{burlet derham}, \cite{furuya porto}, \cite{saeki 2} and \cite{sakuma} for example. Let $f$ be a fold
 map. Then the singular set $S(f)$ and the set of all the fold points whose indices are $i$ (we denote the set by $F_i(f)$) are
 $C^{\infty}$ {\rm (}$n-1${\rm)}-submanifolds of $M$. The restriction map $f {\mid}_{S(f)}$ is a $C^{\infty}$ immersion. \\
\ \ \ A Morse function on a closed manifold is a fold map. A Morse function on a closed manifold which has just two singular points is a special generic map. Conversely, a Morse function which is a special generic map
 on a closed and connected manifold whose dimension is larger than $1$ has just two singular points. \\
\ \ \ In this paper, we sometimes need Morse functions on compact $C^{\infty}$ manifolds possibly with boundaries. We call
 a Morse function on such a manifold {\it good} if it is constant and minimal on the boundary, singular
 points of it are not on the boundary and at any two distinct singular points, the values are distinct. \\
\ \ \ We introduce {\it stable} fold maps.  

\begin{Def}
\label{def:2}
A fold map $f:M \rightarrow N$ is said to be a {\it stable} fold map if the restriction $f {\mid}_{S(f)}$ is a $C^{\infty}$ immersion with normal crossings.
\end{Def}

Note that a stable fold map on a closed $C^{\infty}$ manifold is also a {\it stable} map (for stable maps, see \cite{golubitsky guillemin} for example). \\
\ \ \ Note also that a Morse function on a closed $C^{\infty}$ manifold is a stable fold map if and only if it is good. \\
\ \ \ We also introduce {\it simple} fibers of fold maps and {\it simple} fold maps (see also \cite{saeki} and \cite{sakuma} for example). 

\begin{Def}
\label{def:3}
For a fold map $f$ and $p \in f(S(f))$ , $f^{-1}(p)$ is said to be {\it simple} if each connected component of $f^{-1}(p)$ includes
 at most one singular point of $f$. $f$ is said to be a {\it simple} fold map if for each $p \in f(S(f))$, $f^{-1}(p)$
 is simple.   
\end{Def}

\begin{Ex}
\label{ex:1}
\begin{enumerate}
\item A Morse function on a closed manifold is simple if it is good. 
\item A fold map $f:M \rightarrow {\mathbb{R}}^n$ is simple if $f {\mid}_{S(f)}$ is a $C^{\infty}$ embedding.  
\item Special generic maps are simple.
\end{enumerate}
\end{Ex}

\subsection{Reeb spaces}

We review the {\it Reeb space} of a map. 

\begin{Def}
\label{def:4}
 Let $X$, $Y$ be topological spaces. For $p_1, p_2 \in X$ and for a map $c:X \rightarrow Y$, 
 we define as $p_1 {\sim}_c p_2$ if and only if $p_1$ and $p_2$ are in
 the same connected component of $c^{-1}(p)$ for some $p \in Y$. ${\sim}_{c}$ is an equivalence relation. \\
\ \ \ We denote the quotient space $X/{\sim}_c$ by $W_c$. We call $W_c$ the {\it Reeb space} of $c$.
\end{Def}

 We denote the induced quotient map from $X$ into $W_c$ by $q_c$. We define $\bar{c}:W_c \rightarrow Y$
 so that $c=\bar{c} \circ q_c$. $W_c$ is often homeomorphic to a polyhedron. \\
\ \ \ For example, for a
 good Morse function, the Reeb space
 is a graph. For a simple fold map, the Reeb space is homeomorphic to a polyhedron which is not so complex (see Proposition \ref{prop:2}). \\
\ \ \ For a
 special generic map,
 the Reeb space is homeomorphic to a $C^{\infty}$ manifold. See section 2 of \cite{saeki 2}. See also \cite{burlet derham} and \cite{furuya porto} for example. \\ 
\ \ \ In \cite{kobayashi saeki}, it is proven that for a stable fold map (or stable map) $f:M \rightarrow {\mathbb{R}}^2$ ($m \geq 2$), $W_f$ is homeomorphic to a polyhedron.
 It is known 
that for a stable map $f$, the Reeb space $W_f$ is homeomorphic to a polyhedron. See \cite{hiratuka} for example. \\
\ \ \ The following holds since a stable fold map is stable.

\begin{Prop}
\label{prop:1}
For a stable fold map $f$, the Reeb space $W_f$ is homeomorphic to a polyhedron.
\end{Prop}

The following Proposition \ref{prop:2} is well-known and we omit the proof. See \cite{kobayashi saeki}, \cite{saeki} and
 \cite{saeki suzuoka} for example. In this paper, we often apply the statements of this proposition implicitly. \\
\ \ \ We also note that in this paper, an {\it almost-sphere} of dimension $k$ means a $C^{\infty}$ homotopy sphere given by glueing two standard
 closed discs
 of dimensions $k$ together by a $C^{\infty}$ diffeomorphism between the boundaries. \\
\ \ \ We often use terminologies on (fiber) bundles in this paper (see also \cite{steenrod}). For a topological space $X$, an {\it $X$-bundle} is a bundle
 whose fiber is $X$. A bundle whose structure group is $G$
 is said to be a {\it trivial} bundle if it is equivalent to the product bundle as a bundle whose structure group is $G$. Especially, a trivial bundle whose
 structure group is a subgroup of the homeomorphism group of the fiber is said to be a {\it topologically trivial} bundle. In this paper, a {\it $C^{\infty}$ {\rm (}${\rm PL}${\rm )} bundle} means a bundle
 whose fiber is a $C^{\infty}$ (resp. ${\rm PL}$) manifold
 and whose structure group is a subgroup of the $C^{\infty}$ diffeomorphism group (resp. ${\rm PL}$ homeomorphism group) of the
 fiber. A {\it linear} bundle is
 a $C^{\infty}$ bundle whose fiber is a standard disc or a standard sphere or ${\mathbb{R}}^n (n \geq 1)$ and whose structure group is a subgroup of an orthogonal group.

\begin{Prop}
\label{prop:2}
Let $f:M \rightarrow N$ be a special generic map or a simple fold map or a stable fold map. Then $W_f$ has the structure of a polyhedron and the followings hold.  
\begin{enumerate}
\item $W_f-q_f(S(f))$ is uniquely given the structure of a $C^{\infty}$ manifold such that ${q_f} {\mid}_{M-S(f)}:M-S(f) \rightarrow W_f-q_f(S(f))$ is a $C^{\infty}$ submersion. Furthermore, for any
 compact $C^{\infty}$ submanifold $R$ of dimension $n$
 of any connected component of $W_f-q_f(S(f))$, $R$ is a subpolyhedron of $W_f$ and $q_f {\mid}_{{q_f}^{-1}(R)}:{q_f}^{-1}(R) \rightarrow R$ gives the structure of a $C^{\infty}$ bundle whose fiber is a connected $C^{\infty}$
 manifold of dimension $m-n$.
\item The restriction of $q_f$ to the set $F_0(f)$ of all the definite fold points, is injective.
\item $f$ is simple if and only if ${q_f}{\mid}_{S(f)}:S(f) \rightarrow W_f$ is injective. Special generic maps are simple.     
\item If $f$ is simple, then for any connected component $C$ of $S(f)$, $q_f(C)$ has a small
 regular neighborhood $N(q_f(C))$ in $W_f$
 and ${q_f}^{-1}(N(q_f(C)))$ has the structure of a $C^{\infty}$ bundle
 over $q_f(C)$.   
\item For any connected component $C$ of $F_0(f)$, any small
 regular neighborhood of $q_f(C)$ has the structure of a trivial ${\rm PL}$ $[0,1]$-bundle over $q_f(C)$ and $q_f(C)$ corresponds to
 the $0$-section. We can take
 a small regular neighborhood $N(q_f(C))$ of $q_f(C)$
 and ${q_f}^{-1}(N(q_f(C)))$ has the structure of a linear $D^{m-n+1}$-bundle
 over $q_f(C)$. More precisely, the bundle structure is given as the following; for the connected
 component $C^{\prime}$ of $\partial N(q_f(C))$ satisfying $C^{\prime} \bigcap q_f(F_0(f))=\emptyset$,
 the composition of $q_f {\mid}_{{q_f}^{-1}(N(q_f(C))}:{q_f}^{-1}(N(q_f(C))) \rightarrow N(q_f(C))$ and the projection
 to $q_f(C)$ {\rm (}the subbundle
 corresponding to the $0$-section{\rm )} or ${C}^{\prime}$ {\rm (}the subbundle
 corresponding to the fiber $\{1\} \subset [0,1]${\rm )} .
\item Let $f$ be simple. Let $m-n \geq 1$. If $m-n=1$, then we also assume that $M$ is orientable. \\
\ \ \ Then for any connected component $C$ of the set $F_1(f)$ of all
 the fold points of indice $1$, such that for any connected component $R$ of $W_f-q_f(S(f))$
 whose closure $\overline{R}$ includes $q_f(C)$, ${q_f}^{-1}(p)$ is an almost-sphere for $p \in R$, any small
 regular neighborhood of $q_f(C)$ has the structure of a $K$-bundle over $q_f(C)$ where $K:= \{r \exp (2\pi i \theta) \in \mathbb{C} \mid 0 \leq r \leq 1, \theta = 0, \frac{1}{3}, \frac{2}{3} \}$ , where the structure group consists of just two elements fixing the point $1 \in K$ and where $q_f(C)$ corresponds to the $0$-section.
 We can take a small regular neighborhood $N(q_f(C))$ of $q_f(C)$
 and ${q_f}^{-1}(N(q_f(C)))$ has the structure
 of a $C^{\infty}$ bundle over $q_f(C)$ whose fiber is ${\rm PL}$ homeomorphic to $S^{m-n+1}$ with the interior
 of a union of disjoint three standard
 closed {\rm (}$m-n+1${\rm )}-discs removed. More precisely, the bundle structure is given by the composition of $q_f {\mid}_{{q_f}^{-1}(N(q_f(C)))}:{q_f}^{-1}(N(q_f(C))) \rightarrow N(q_f(C))$ and
 the projection to $q_f(C)$ {\rm (}the subbundle
 corresponding to the $0$-section{\rm )} or the projection to the subbundle corresponding to the fiber $\{1\} \subset K$.     
\end{enumerate}  
\end{Prop}

\section{Notes on round fold maps}
\label{sec:3}

In this section, we review round fold maps. See also \cite{kitazawa 2}.

\subsection{Terms on round fold maps}

\begin{Def}[round fold maps (\cite{kitazawa 2})]
\label{def:5}
$f:M \rightarrow {\mathbb{R}}^n$ ($m \geq n \geq 2$) is said to be a {\it round} fold map if $f$ is $C^{\infty}$ equivalent to
 a fold map $f_0:M_0 \rightarrow {\mathbb{R}}^n$ on a closed $C^{\infty}$ manifold $M_0$ such that the followings hold.

\begin{enumerate}
\item The singular set $S(f_0)$ is a disjoint union of standard spheres of dimensions $n-1$ and consists of $l \in \mathbb{N}$ connected components.
\item The restriction map $f_0 {\mid}_{S(f_0)}$ is a $C^{\infty}$ embedding.
\item Let ${D^n}_r:=\{(x_1,\cdots,x_n) \in {\mathbb{R}}^n \mid {\sum}_{k=1}^{n}{x_k}^2 \leq r \}$. Then $f_0(S(f_0))={\sqcup}_{k=1}^{l} \partial {D^n}_k$.  
\end{enumerate}

We call $f_0$ a {\it normal form} of $f$. We call a ray $L$ from $0 \in {\mathbb{R}}^n$ an {\it axis} of $f_0$ and
 ${D^n}_{\frac{1}{2}}$ the {\it proper core} of $f_0$. Suppose that for a round fold map $f$, its normal form $f_0$ and $C^{\infty}$ diffeomorphisms
 $\Phi:M \rightarrow M_0$ and $\phi:{\mathbb{R}}^n \rightarrow {\mathbb{R}}^n$, $\phi \circ f=f_0 \circ \Phi$. Then
 for an axis $L$ of $f_0$, we also call ${\phi}^{-1}(L)$ an {\it axis} of $f$ and for the proper core ${D^n}_{\frac{1}{2}}$ of $f_0$, we
 also call ${\phi}^{-1}({D^n}_{\frac{1}{2}})$ a {\it proper core} of $f$. 
\end{Def}

\begin{figure}
\begin{center}
\includegraphics[width=50mm]{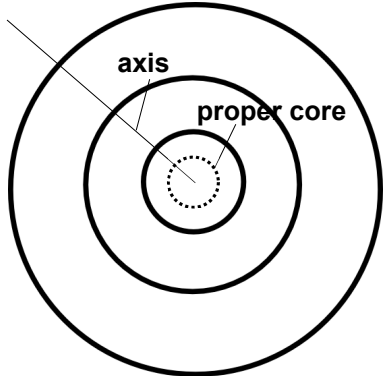}
\end{center}
\caption{An axis and a proper core of a round fold map}
\label{fig:1}
\end{figure}

We introduce another definition.

\begin{Def}[round fold maps (\cite{kitazawa 2})]
\label{def:6}
Assume that $f:M \rightarrow {\mathbb{R}}^n$ is a fold map and that $m \geq n \geq 2$.
 $f$ is said to be a {\it round} fold map if the followings hold.

\begin{enumerate}
\item The singular set $S(f)$ is a disjoint union of standard spheres of dimensions $n-1$.
\item The restriction map $f {\mid}_{S(f)}$ is a $C^{\infty}$ embedding. 
\item We denote by $\{U_0\} \sqcup \{U_{\infty}\} \sqcup {\{U_{\lambda}\}}_{\lambda \in \Lambda}$ ($\Lambda$ may be empty.)  the set of all the connected components
 of ${\mathbb{R}}^n-f(S(f))$. The followings hold. 
\begin{enumerate}
\item The closure $\overline{U_0}$ is $C^{\infty}$ diffeomorphic to $D^n$. 
\item The closure $\overline{U_{\infty}}$ is $C^{\infty}$ diffeomorphic to $S^{n-1} \times [0,+\infty)$.
\item The closure $\overline{U_{\lambda}}$ is $C^{\infty}$ diffeomorphic to $S^{n-1} \times [0,1]$. 
\end{enumerate}
\end{enumerate}
\end{Def}

\begin{figure}
\begin{center}
\includegraphics[width=50mm]{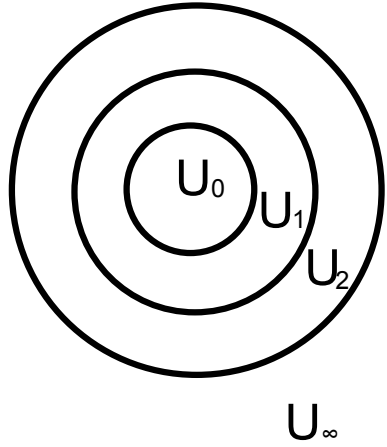}
\end{center}
\caption{The image of a round fold map (Definition \ref{def:6})}
\label{fig:2}
\end{figure}

In \cite{kitazawa 2}, the following has been shown.

\begin{Prop}[\cite{kitazawa 2}]
\label{prop:3}
Two definitions of a round fold map are equivalent.
\end{Prop}

Let $f:M \rightarrow {\mathbb{R}}^n$ be a round fold map, $L$ be its axis and $P$ be its proper core. \\
\ \ \ We set $E:=f^{-1}(P)$ and $E^{\prime}:=M-f^{-1}({\rm Int} P)$. We put
 $F:=f^{-1}(p)$
 for $p \in \partial P$. \\
\ \ \ We put $P^{(1)}:=P$, $P^{(2)}:={\mathbb{R}}^n-{\rm Int} P$,
 $f_1:=f {\mid}_{E}:E \rightarrow P^{(1)}$ and $f_2:=f {\mid}_{{E}^{\prime}}:{E}^{\prime} \rightarrow P^{(2)}$. \\
\ \ \ $f_1$ gives the structure of a trivial $C^{\infty}$ bundle over $P^{(1)}$, ${f_1} {\mid}_{\partial E}:\partial E \rightarrow \partial P^{(1)}$ gives
 the structure of a trivial $C^{\infty}$ bundle
 over $\partial P^{(1)}$ and $f_2 {\mid}_{\partial {E}^{\prime}}:\partial {E}^{\prime} \rightarrow \partial P^{(2)}$
 gives the structure of a trivial $C^{\infty}$ bundle over $\partial P^{(2)}$ if $P^{(1)} \bigcap f(M) \neq \emptyset$. \\
\ \ \ We can give ${E}^{\prime}$ and $q_f({E}^{\prime})$ the structures of bundles over $\partial {P}^{(2)}$ as follows. \\ 
\ \ \ Let $f$ be a normal form of a round fold map. Since for ${\pi}_P(x):=\frac{1}{2} \frac{x}{|x|}$ ($x \in P^{(2)}$), ${\pi}_P \circ f {\mid}_{{E}^{\prime}}$ is
 a proper $C^{\infty}$ submersion, this map gives ${E}^{\prime}$ the structure of
 a $C^{\infty}$ $f^{-1}(L)$-bundle over $\partial P^{(2)}$ (apply Ehresmann's fibration
 theorem; see also \cite{ehresmann}). \\
\ \ \ $q_f({E}^{\prime})$ has the structure of a bundle over $\partial P^{(2)}$. We consider the following map. \\
\ \ \ Let $p \in q_f({E}^{\prime})$. For $p_1,p_2 \in {q_f}^{-1}(p)$, ${\pi}_P \circ f(p_1)={\pi}_P \circ f(p_2)$. We correspond
 ${\pi}_P \circ f(p_1)={\pi}_P \circ f(p_2)$ to $p$. This map from $q_f({E}^{\prime})$ into $\partial P^{(2)}$
 gives the structure of a ${\bar{f}}^{-1}(L)$-bundle: we only notice that ${\pi}_P \circ f {\mid}_{{E}^{\prime}}$ gives the structure of a $C^{\infty}$ bundle
 and that $q_f({E}^{\prime})$ is the quotient space of ${E}^{\prime}$ by ${\sim}_f$. \\
\ \ \ For a round fold map $f$ which is not a normal form, we can consider similar structures of bundles. \\ 

\begin{Def}[\cite{kitazawa 2}]
\label{def:7}
Let $f:M \rightarrow {\mathbb{R}}^n$ be a round fold map. 
If the natural projection from the total space of the bundle ${E}^{\prime}$ onto the base space $\partial P^{(2)}$ ($\partial P^{(2)}$ is $C^{\infty}$ diffeomorphic to $S^{n-1}$) gives
 the structure of a topologically ($C^{\infty}$) trivial bundle, then
 $f$ is said to be {\it topologically {\rm (}{\rm resp.} $C^{\infty}${\rm )} trivial}. \\
\end{Def}

Related with topologically or $C^{\infty}$ trivial round fold maps and differential topology of manifolds admitting
 such maps, we introduce Theorem 1 of \cite{kitazawa 4} in the following. 

\begin{Prop}[\cite{kitazawa 4}]
\label{prop:4}
Let $M$ be a closed $C^{\infty}$ manifold of dimension $m$. Let $n \in \mathbb{N}$ and $m \geq n \geq 2$. 
\begin{enumerate}
\item
\label{prop:4.1}
 Let $M$ have the structure of a $C^{\infty}$ bundle over $S^n$ whose fiber is a closed $C^{\infty}$ manifold $F (\neq \emptyset)$. Then there exists
 a $C^{\infty}$ trivial round fold map $f:M \rightarrow {\mathbb{R}}^n$ such that the fiber of a point in a proper core of $f$ is $C^{\infty}$ diffeomorphic
 to a disjoint union of two copies of $F$ and that for an axis $L$ of $f$, $f^{-1}(L)$ is $C^{\infty}$ diffeomorphic to $F \times [0,1]$. 
\item
\label{prop:4.2}
 Suppose that a topologically {\rm (}$C^{\infty}${\rm )} trivial
 round fold map $f:M \rightarrow {\mathbb{R}}^n$ exists and that for an axis $L$ of $f$ and a closed $C^{\infty}$ manifold $F$ of dimension $m-n$, $f^{-1}(L)$ is $C^{\infty}$ diffeomorphic
 to $F \times [0,1]$. Then $M$ has the structure of a {\rm (}resp. $C^{\infty}${\rm )} $F$-bundle over $S^n$. 
\end{enumerate}
\end{Prop}

\ref{prop:4.1} of Proposition \ref{prop:4} has been again shown as the following as Theorem 3 of \cite{kitazawa 5}.

\begin{Prop}[\cite{kitazawa 5}]
\label{prop:5}
Let $F \neq \emptyset$ be a closed and connected $C^{\infty}$ manifold. Let $M$ be a closed $C^{\infty}$ manifold of dimension $m$ having the structure of a $C^{\infty}$ $F$-bundle
 over $S^n$ {\rm (}$n \geq 2${\rm )}. Then $M$ admits a round fold map $f:M \rightarrow {\mathbb{R}}^n$ satisfying the followings.
\begin{enumerate}
\item $f$ is $C^{\infty}$ trivial.
\item For an axis $L$ of $f$, $f^{-1}(L)$ is $C^{\infty}$ diffeomorphic to $F \times [0,1]$.
\item Two connected components of the fiber of a point in a proper core of $f$ is regarded as fibers of the $C^{\infty}$ $F$-bundle
 over $S^n$.
\item $f(M)$ is $C^{\infty}$ diffeomorphic to $D^n$ and for the connected component $C:=\partial f(M)$, the $C^{\infty}$ embedding
 of $f^{-1}(C)$ into $M$ is $C^{\infty}$ isotopic to a trivial $C^{\infty}$ embedding into ${\rm Int} D^m {\rm (}\subset D^m \subset M_2{\rm )}$.
\end{enumerate}
\end{Prop}
 	
\section{Round fold maps on bundles over manifolds admitting round fold maps}
\label{sec:4}

In this section, we extend or generalize statements of Proposition \ref{prop:4} and Proposition \ref{prop:5}. More precisely, we construct round fold maps on product bundles over manifolds admitting
 {\it locally $C^{\infty}$ trivial} round fold maps and on $C^{\infty}$ bundles over (exotic) spheres in new manners. \\
\ \ \ First
 we introduce a {\it locally $C^{\infty}$ trivial} round fold map.

\begin{Def}
\label{def:8}
Let $f:M \rightarrow {\mathbb{R}}^n$ be a round fold map. Assume that for any connected component $C$ of $f(S(f))$ and a small $C^{\infty}$ closed tubular neighborhood $N(C)$ of $C$ such that $\partial N(C)$ is the disjoint
 union of two connected components $C_1$ and $C_2$, $f^{-1}(N(C))$ has the structure
 of a trivial $C^{\infty}$ bundle over $C_1$ {\rm (}$C_2${\rm )} and $f {\mid}_{f^{-1}(C_1)}:f^{-1}(C_1) \rightarrow C_1$ {\rm (}resp. $f {\mid}_{f^{-1}(C_2)}:f^{-1}(C_2) \rightarrow C_2${\rm )} gives the structure
 of a subbundle of the bundle $f^{-1}(N(C))$. Then $f$ is said to be {\it locally $C^{\infty}$ trivial}. We call a fiber $F_C$ of the bundle $f^{-1}(N(C))$ the {\it normal fiber of $C$ corresponding to the bundle $f^{-1}(N(C))$}. Assume that $C_1$ is in the bounded connected
 component of ${\mathbb{R}}^n-C_2$ and we denote a fiber of the subbundle $f^{-1}(C_1)$ by ${\partial}_0 F_C$.
\end{Def}

Note that a round fold map is locally $C^{\infty}$ trivial if it is $C^{\infty}$ trivial. \\
\ \ \ We introduce results shown in \cite{kitazawa 4} and \cite{kitazawa 5}. 
Note that a locally $C^{\infty}$ trivial map is first defined in this paper, although maps satisfying
 the definition of a locally $C^{\infty}$ trivial map appeared in these papers. \\
\ \ \ We also review a {\it round special generic group} defined in Definition 11 of \cite{kitazawa 4}. Let ${\Theta}_{k_1}$ be the h-cobordism group of $C^{\infty}$ oriented homotopy spheres of dimensions $k_1 \geq 2$. It follows
 easily that the set
 of all the classes of ${\Theta}_{k_1}$ consisting of spheres admitting round fold maps with connected singular sets
 into ${\mathbb{R}}^{k_2}$ ($k_1 \geq k_2 \geq 2$) is a subgroup of ${\Theta}_{k_1}$. In fact we only consider the connected
 sum of given two round fold maps with connected singular sets (for the connected sum of given two special generic maps into Euclidean spaces,
 see section 5 of \cite{saeki} for example). We denote the subgroup by ${\Theta}_{(k_1,k_2)} \subset {\Theta}_{k_1}$ and call it the {\it $(k_1,k_2)$ round special generic group}.

\begin{Prop}[\cite{kitazawa 4}, \cite{kitazawa 5}]
\label{prop:6}
\begin{enumerate}
\item
\label{prop:6.1}
 Any connected sum of $l \in \mathbb{N}$ closed $C^{\infty}$ manifolds of dimensions $m$ having the structures of $C^{\infty}$ $S^{m-n}$-bundles over $S^n$ admits a locally $C^{\infty}$ trivial round fold map $f$ into ${\mathbb{R}}^n$ satisfying the followings where $n \geq 2$ and $m-n \geq 1$.
\begin{enumerate}
\item All the regular fibers of $f$ are disjoint unions of finite copies of $S^{m-n}$.
\item The number of connected components of $S(f)$ and the number of connected components of the fiber of a point in a proper core of $f$ are $l$. 
\item All the connected components of the fiber of a point in a proper core of $f$ are regarded as fibers of the $C^{\infty}$ $S^{m-n}$-bundles
 over $S^n$ and a fiber of any $C^{\infty}$ $S^{m-n}$-bundle over $S^n$ appeared in the connected sum is regarded as a connected component of the
 fiber of a point in a proper core of $f$.
\end{enumerate}
\item
\label{prop:6.2}
 Let $M$ be a closed and connected $C^{\infty}$ oriented manifold of dimension $m$. Let $n \in \mathbb{N}$, $n \geq 2$ and $m \geq 2n$. Then
 the followings are equivalent.
\begin{enumerate}
\item 
 A locally $C^{\infty}$ trivial round fold map $f:M \rightarrow {\mathbb{R}}^n$ exists and the followings hold.
\begin{enumerate}
\item Regular fibers of $f$ are disjoint
 unions of standard spheres.
\item The number of connected components of the fiber of a point in a proper core of $f$ equals
 the number of connected components of $S(f)$. 
\end{enumerate}
\item
$M$ is a connected sum of a finite number of $C^{\infty}$ oriented manifolds of dimensions $m$ having the structures of $C^{\infty}$ $S^{m-n}$-bundles over $S^n$ and an oriented manifold in a class of ${\Theta}_{(m,n)}$.
\end{enumerate}
\end{enumerate}
\end{Prop}

\begin{Prop}
\label{prop:7}
Let $M$ be a closed $C^{\infty}$ manifold of dimension $m$.
Assume that there exists a {\rm (}locally{\rm )} $C^{\infty}$ trivial round
 fold map $f:M \rightarrow {\mathbb{R}}^n$ {\rm (}$m \geq n \geq 2${\rm )}. Then for any closed $C^{\infty}$ manifold $F$, we
 can construct a {\rm (}resp. locally{\rm )} $C^{\infty}$ trivial round fold map ${f}^{\prime}:M \times F \rightarrow {\mathbb{R}}^n$
\end{Prop}
\begin{proof}
We prove the case where $f$ is locally $C^{\infty}$ trivial. \\
\ \ \ We may assume that $f:M \rightarrow {\mathbb{R}}^n$ is a normal form. Let $S(f)$ consist of $l$ connected
 components. Set $P_k:={D^n}_{k+\frac{1}{2}}-{\rm Int} {D^n}_{k-\frac{1}{2}}$ for an integer $1 \leq k \leq l$. Then we may assume
 that ${f}^{-1}(P_k)$ has the structure of a trivial $C^{\infty}$ bundle
 over $\partial {D^n}_{k-\frac{1}{2}}$ ($\partial {D^n}_{k+\frac{1}{2}}$) with a fiber $C^{\infty}$ diffeomorphic to $E_k$ such
 that $f {\mid}_{{f^{-1}}(\partial {D^n}_{k-\frac{1}{2}})}$ (resp. $f {\mid}_{f^{-1}(\partial {D^n}_{k+\frac{1}{2}})}$) gives
 the structure of subbundle (we denote fibers of two subbundles by ${E_k}^{1} \subset E_k$ and ${E_k}^2 \subset E_k$, respectively). For a $C^{\infty}$ diffeomorphism ${\phi}_k$ regarded as a
 bundle isomorphism between the two trivial $C^{\infty}$ bundles (over standard spheres) inducing the identification between the base spaces, $M$ is $C^{\infty}$ diffeomorphic to
 $$(\cdots ((f^{-1}({D^n}_{\frac{1}{2}})) {\bigcup}_{{\phi}_1} f^{-1}(P_1)) \cdots) {\bigcup}_{{\phi}_l} f^{-1}(P_l)$$
 and $M \times F$ is
 regarded as
 $$(\cdots ((f^{-1}({D^{n}}_{\frac{1}{2}}) \times F) {\bigcup}_{{{\phi}_1} \times {\rm id}_F} (f^{-1}(P_1) \times F)) \cdots) {\bigcup}_{{\phi}_l \times {\rm id}_F} ({f^{-1}}(P_l) \times F)$$

We construct a $C^{\infty}$ map on $f^{-1}(P_k) \times F$. This manifold has the structure of a
 trivial $C^{\infty}$ bundle over $\partial {D^n}_{k-\frac{1}{2}}$ (and $\partial {D^n}_{k+\frac{1}{2}}$) with
 a fiber $C^{\infty}$ diffeomorphic to $E_k \times F$. On $E_k \times F$ there exists a Morse function $\tilde{f_k}$ and the followings hold.

\begin{enumerate}
\item $\tilde{f_k}$ is constant and minimal on $\partial {E_k}^1 \times F$. (We denote the minimum by $a$.) 
\item $\tilde{f_k}$ is constant and maximal on $\partial {E_k}^2 \times F$. (We denote the maximum by $b$.)
\item Singular points of $\tilde{f_k}$ are in the interior of $E_k \times F$ and at any two distinct singular points, the values are different.
\end{enumerate}

We obtain ${\rm id}_{S^{n-1}} \times \tilde{f_k}:S^{n-1} \times E_k \times F \rightarrow S^{n-1} \times [a,b]$.
 By glueing $f {\mid}_{f^{-1}({D^n}_{\frac{1}{2}})} \times {\rm id}_F$ and the family
 $\{{\rm id}_{S^{n-1}} \times \tilde{f_k}\}$ by using the family $\{{{\phi}_k} \times {\rm id}_{S^{n-1}} \}$ and the family of identifications (on the target), we obtain
 a new round fold map ${f}^{\prime}:M \times F \rightarrow {\mathbb{R}}^n$. 
\end{proof}

\begin{Cor}
\label{cor:1}
For any closed $C^{\infty}$ manifold $F$, any $C^{\infty}$ oriented homotopy sphere $\Sigma$ of any class in ${\Theta}_{m,n}$ {\rm (}$m \geq n \geq 2${\rm )}, $\Sigma \times F$ admits a $C^{\infty}$ trivial round fold map into ${\mathbb{R}}^n$ such that the fiber of a point in
 a proper core is $S^{m-n} \times F$. 
\end{Cor}
\begin{proof}
We only note that a round fold map with connected singular set from $\Sigma$ into ${\mathbb{R}}^n$ is $C^{\infty}$ trivial since the natural projection in Definition \ref{def:7} gives the
structure of a linear bundle and the projection gives the structure of a
trivial $C^{\infty}$ bundle on the boundary (see also Example 3 (1) of \cite{kitazawa 2}, for example). We only apply the method in the proof
 of Proposition \ref{prop:7} to the round fold map.
\end{proof}

\begin{Rem}
Any almost-sphere $\Sigma$ of dimension larger than $1$ admits a round fold map into the plane with connected singular set and the map is $C^{\infty}$ trivial. We can apply the method in the proof on the map to construct a $C^{\infty}$ trivial round fold map from $\Sigma \times F$ into ${\mathbb{R}}^2$ such that the fiber of a point in
 a proper core is $C^{\infty}$ diffeomorphic to $S^{m-2} \times F$.
\end{Rem}

\begin{Cor}
\label{cor:2}
Let $M$ be a closed $C^{\infty}$ manifold of
 dimension $m$. Let $F_1$ be a closed $C^{\infty}$ manifold of dimension $m-n$ {\rm (}$m \geq n \geq 2${\rm )}. Let $M$ have
 the structure of $C^{\infty}$ $F_1$-bundle over $S^n$. Then for any closed $C^{\infty}$ manifold $F_2$, there exists
 a $C^{\infty}$ trivial round fold map ${f}^{\prime}:M \times F_2 \rightarrow {\mathbb{R}}^n$ such that
 the fiber of a point in a proper core is $C^{\infty}$ diffeomorphic
 to a disjoint union of two copies of $F_1 \times F_2$. 
\end{Cor}
\begin{proof}
We only apply the construction in the proof of Proposition \ref{prop:7} to a round fold map from $M$ into ${\mathbb{R}}^n$ obtained in Proposition \ref{prop:4} or Proposition \ref{prop:5}. 
\end{proof}

\begin{Thm}
\label{thm:1}
For any $C^{\infty}$ oriented homotopy sphere $\Sigma$ in a class of ${\Theta}_{m,n}$ where $m \geq n \geq 2$ and any closed $C^{\infty}$ manifold $M$ having the structure
 of a $C^{\infty}$ $F$-bundle
 over $\Sigma$ admits a $C^{\infty}$ trivial round
 fold map into ${\mathbb{R}}^n$.    
\end{Thm}
\begin{proof}
$\Sigma$ admits a round fold map $f$ whose singular set is connected into ${\mathbb{R}}^n$ and the map $f$ is $C^{\infty}$ trivial. Let $P$ be a proper
 core of $f$ and $Q$ be a small $C^{\infty}$ closed tubular neighborhood of the singular value set $f(S(f))$. We may assume
 that $\partial P \subset \partial Q$. For a $C^{\infty}$ diffeomorphism $\phi$ which is regarded as an isomorphism between the two trivial
 $C^{\infty}$ bundles $\partial P \times S^{m-n}$ over $\partial P$ and $\partial P \times \partial D^{m-n+1}$ over $\partial P$ inducing the
 identification between the base spaces, $\Sigma$ is regarded as $(P \times S^{m-n}) {\bigcup}_{\phi} (\partial P \times D^{m-n+1})$. If we restrict
 the base space of any $C^{\infty}$ $F$-bundle over $\Sigma$ to $P \times S^{m-n} \subset \Sigma$ and $\partial P \times D^{m-n+1} \subset \Sigma$, then
 the resulting
 $C^{\infty}$ bundles are trivial $C^{\infty}$ bundles. So for every $C^{\infty}$ $F$-bundle over $\Sigma$, there exists a $C^{\infty}$ diffeomorphism $\Phi$ regarded as a bundle isomorphism
 between the two trivial $C^{\infty}$ $F$-bundles
 $\partial P \times S^{m-n} \times F$ over $\partial P \times S^{m-n}$ and $\partial P \times \partial D^{m-n+1} \times F$ over
 $\partial P \times \partial D^{m-n+1}$ inducing $\phi$ and the $C^{\infty}$ $F$-bundle over $\Sigma$ is regarded as $(P \times S^{m-n} \times F) {\bigcup}_{\Phi} (\partial P \times D^{m-n+1} \times F)$. \\
\ \ \ We have a good Morse function $\tilde{f}:D^{m-n+1} \times F \rightarrow [a,+\infty)$ where $a$ is the minimum.
 By using $\Phi$ and the identification map on the target, we can glue the canonical projection $p:P \times S^{m-n} \times F \rightarrow P$ and $\tilde{f} \times {\rm id}_{\partial P}$ to obtain a desired round fold map $f$. This completes the proof.    
\end{proof}

\begin{Rem}
 Theorem \ref{thm:1} has been shown as a part of Theorem 1 of \cite{kitazawa 4} or Proposition \ref{prop:4} of the present paper for the case
 where $\Sigma=S^n$ (in section 4 of \cite{kitazawa}, it has been shown for the case where $\Sigma=S^n$ and where $F$ is an almost-sphere). The proof in the present paper
 can be regarded as an extension of the original proofs. 
\end{Rem}

Now we introduce terms and results of \cite{kitazawa 3}. \\
\ \ \ Before that, we define a {\it branched} point of a polyhedron. For a polyhedron of dimension $k \geq 1$, 
 a {\it branched} point means a point such that every open
 neighborhood of the point is not homeomorphic to any open set of ${\mathbb{R}}^k$ or ${{\mathbb{R}}^k}_{+}:=\{(x_1,\cdots,x_k) \in {\mathbb{R}}^k \mid x_k \geq 0\}$. If a polyhedron $X$ of dimension
 $k$ does not have branched points, then it is a manifold having a triangulation and we can define the {\it interior} ${\rm Int} X$ and the {\it boundary} $\partial X$.   \\

\begin{Def}
\label{def:9}
Let $f:M \rightarrow N$ be a simple fold map and $m-n \geq 2$.

\begin{enumerate}
\item Let $C$ be a connected component of $q_f(S(f))-q_f(F_0(f))$ consisting of non-branched points. Assume that there exists a small regular
 neighborhood $N(C)$ of $C$ in $W_f$ having the
 structure of a trivial ${\rm PL}$ $[-1,1]$-bundle over $q_f(C)$ ($q_f(C)$ corresponds to the $0$-section) and
 that the composition of ${q_f} {\mid}_{{q_f}^{-1}(N(C))}:{q_f}^{-1}(N(C)) \rightarrow N(C)$ and the projection to $q_f(C)$ gives ${q_f}^{-1}(N(C))$ the structure
 of a $C^{\infty}$ bundle over $C$ whose fiber is ${\rm PL}$ homeomorphic
 to $D^{k+1} \times S^{m-n-k}-{\rm Int} D^{m-n+1}$ for an integer $1 \leq k < m-n$. 
Then $C$ is said to be a {\it $k$ S-locus}.
\item Let $R$ be a connected component of $W_f-q_f(S(f))$. \\
\ \ \ If for any compact subset $P$ in $R$,
 ${q_f} {\mid}_{{q_f}^{-1}(P)}:{q_f}^{-1}(P) \rightarrow P$ gives the structure of a bundle over $P$ whose fiber is ${\rm PL}$ homeomorphic to $S^k \times S^{m-n-k}$ for an integer $1 \leq k \leq [\frac{m-n}{2}]$, 
then $R$ is said to be a {\it $k$ S-region} or an {\it $m-n-k$ S-region} (note that $k=m-n-k$ if $m-n$ is even and $k=\frac{m-n}{2}$). \\
\ \ \ If for $p \in R$, ${q_f}^{-1}(p)$ is an almost-sphere, then $R$ is said to be an {\it AS-region}.   
\end{enumerate}
\end{Def}

\ \ \ For example, on $S^{k_1} \times S^{k_2}$ ($k_1, k_2 \geq 2$), there exists
 a good Morse function $f:S^{k_1} \times S^{k_2} \rightarrow \mathbb{R}$ with four singular points such that $W_f$ contains two $k_1-1$ ($k_2-1$) $S$-loci. \\
\ \ \ Let $f:M \rightarrow N$ be a simple fold map and $p \in S(f)$. Let $q_f(p)$ be in a $k$ S-locus. If $k+1 \leq [\frac{m-n+1}{2}]$, then $p \in F_{k+1}(f)$ and if $k+1>[\frac{m-n+1}{2}]$, then $p \in F_{m-n-k}(f)$. 

\begin{Def}
\label{def:10}
If a simple fold map $f:M \rightarrow N$ ($m-n \geq 2$) satisfies the followings, then $f$ is said to be a {\it normal simple fold map with regular fibers of two spheres}.
\begin{enumerate}
\item Any connected component of $q_f(S(f))$ in $W_f-q_f(F_0(f))$ consisting of non-branched points is an S-locus . 
\item For any connected component $R$ of $W_f-q_f(S(f))$ such that the intersection of the closure $\overline{R}$ and the set of all
 the branched points of $W_f$ is non-empty, the fiber of each point in $R$ is an almost-sphere.
\end{enumerate}
\end{Def}

For example, on $S^{k_1} \times S^{k_2}$ ($k_1,k_2 \in \mathbb{N}$), there exists
 a good Morse function $f:S^{k_1} \times S^{k_2} \rightarrow \mathbb{R}$ with four singular points which is also a normal simple fold map with regular fibers of two spheres. 

\begin{Def}
\label{def:11}
Let $f$ be a normal simple fold map with regular fibers of two spheres.
\begin{enumerate}
\item If for any $k$ S-locus $C$, the connected component $R$ of
 $W_f-q_f(S(f))$ which is not an AS-region such that $C$ is in the boundary of the closure $\bar{R}$, is
 a $k$ S-region and the boundary of the closure is a disjoint union of $k$ S-loci, then we
 say that {\it $f$ decomposes into S-systems}. \\
\item Assume that $f$ decomposes into S-systems. \\ 
\ \ \ We say that {\it $f$ has a family of S-identifications} if there exists a family of small regular neighborhoods $\{N(C_{\lambda})\}_{\lambda \in \Lambda}$ of
 all the S-loci $\{C_{\lambda}\}_{\lambda \in \Lambda}$ such that the followings hold.
\begin{enumerate}
\item (The bundle structure of ${q_f}^{-1}(N(C_{\lambda}))$ over $C_{\lambda}$) \\
 Let $C_{\lambda}$ be a $k$ S-locus. ${q_f}^{-1}(N(C_{\lambda}))$
 has the structure of a $C^{\infty}$ bundle over $C_{\lambda}$ with the fiber ${\rm PL}$ homeomorphic
 to $D^{k+1} \times S^{m-n-k}-{\rm Int} D^{m-n+1} \subset D^{k+1} \times S^{m-n-k}$ for an integer $1 \leq k < m-n$ such that for the connected component
 $C$ of $\partial N(C_{\lambda})$ in a $k$ S-region, ${q_f}^{-1}(C)$
 has the structure of a subbundle of the bundle ${q_f}^{-1}(N(C_{\lambda}))$ over $C_{\lambda}$ with the fiber ${\rm PL}$ homeomorphic to $\partial D^{k+1} \times S^{m-n-k}$
 and that the structure group of ${q_f}^{-1}(C)$ is a subgroup of
 the $C^{\infty}$ diffeomorphism group of the fiber, consisting of some $C^{\infty}$ diffeomorphisms regarded as bundle isomorphisms on the
 trivial ${\rm PL}$ bundle $\partial D^{k+1} \times S^{m-n-k}$ over $\partial D^{k+1}$ inducing ${\rm PL}$ homeomorphisms on
 the base space $\partial D^{k+1}=S^{k}$ ($1 \leq k < m-n$).
\item (The bundle structure of ${q_f}^{-1}(R)$ over $R$ for a connected component $R$ of $W_f-{\rm Int} {\sqcup}_{\lambda \in \Lambda} N(C_{\lambda})$ in an S-region) \\
For any connected component $R$ of $W_f-{\rm Int} {\sqcup}_{\lambda \in \Lambda} N(C_{\lambda})$ in a $k$ S-region such that the closure of the $k$ S-region is bounded
 by a disjoint union of $k$ S-loci, ${q_f}^{-1}(R)$
 has the structure of a $C^{\infty}$ bundle over $R$ with the fiber ${\rm PL}$ homeomorphic to $S^k \times S^{m-n-k}$ such that the structure group is a subgroup
 of the $C^{\infty}$ diffeomorphism group of the fiber, consisting of some $C^{\infty}$ diffeomorphisms regarded as bundle isomorphisms on the
 trivial ${\rm PL}$ bundle $S^k \times S^{m-n-k}$ over $S^k$ inducing ${\rm PL}$ homeomorphisms on
 the base space $S^k$. 
\item (Identifications of pieces of the source manifold on their boundaries) \\
 For any connected component $R$ of $W_f-{\rm Int} {\sqcup}_{\lambda \in \Lambda} N(C_{\lambda})$ in a $k$ S-region such that the closure of the $k$ S-region is bounded
 by a disjoint union of $k$ S-loci, the identification map of the restriction of a bundle ${q_f}^{-1}(R)$ over $R$ satisfying (b) of this definition to
 any connected component
 $C$ of the boundary $\partial R$ and the subbundle ${q_f}^{-1}(C)$ of a bundle ${q_f}^{-1}(N(C_{\lambda}))$ ($C \subset \partial N(C_{\lambda})$)
 over $C$ satisfying (a) of this definition is a bundle map; the structure groups are a subgroup of
 the $C^{\infty}$ diffeomorphism group of the fiber, consisting of some $C^{\infty}$ diffeomorphisms regarded as bundle isomorphisms on the
 trivial ${\rm PL}$ bundle $\partial D^{k+1} \times S^{m-n-k}$ over $\partial D^{k+1}$ inducing ${\rm PL}$ homeomorphisms on
 the base space $\partial D^{k+1}=S^k$. (Here the base space of the fiber of the bundle satisfying (a) means $\partial D^{k+1}$ and
 the base space of the fiber of the bundle satisfying (b) means $S^k$. We should be careful in the case where $k=m-n-k$.) 
\end{enumerate}
\end{enumerate}
\end{Def}

As results, the following has been shown.

\begin{Prop}[\cite{kitazawa 3}]
\label{prop:8}
 Let $M$ be a closed $C^{\infty}$ manifold of dimension $m$, $N$ be a $C^{\infty}$ manifold of dimension $n$ without
 boundary and $f:M \rightarrow N$ be a simple fold map. Let $n \geq 2$ and $m-n \geq 2$. \\
\ \ \ We assume that $f$ has a family of S-identifications.


Then there exist a compact ${\rm PL}$ manifold $W$ of dimension $m+1$ such that $\partial W = M$, a polyhedron $V$
 and continuous maps $r: W \rightarrow V$ and $s:V \rightarrow W_f$ and the followings hold.
\begin{enumerate}
\item There exist a triangulation of $W$, a triangulation of $V$ and a triangulation of $W_f$ such that $r$ is a simplicial map and that $s$
 is a simplicial map and the followings hold.
\begin{enumerate}
\item For each $p \in V$, $r^{-1}(p)$ collapses to a point and $r$ is a homotopy equivalence.
\item If $p$ is in the closure of an AS-region in $W_f$, then $s^{-1}(p)$ is a point. If $p$ is in an AS-region in $W_f$ and $q \in s^{-1}(p)$, then
 $r^{-1}(q)$ is ${\rm PL}$ homeomorphic to $D^{m-n+1}$. 
\item If $p$ is in a $k$ S-region in $W_f$ whose closure is bounded by a disjoint union of $k$ S-loci, then $s^{-1}(p)$ is ${\rm PL}$ homeomorphic to $S^k$. If
 $p$ is in a $k$ S-region in $W_f$ whose closure is bounded by a disjoint union of $k$ S-loci and $q \in s^{-1}(p)$, then $r^{-1}(q)$ is ${\rm PL}$ homeomorphic
 to $D^{m-n-k+1}$.
\end{enumerate}
\item $W$ collapses to a subpolyhedron ${V}^{\prime}$ such that $r {\mid}_{{V}^{\prime}}:{V}^{\prime} \rightarrow V$ is a ${\rm PL}$ homeomorphism.
\end{enumerate}

If $M$ is orientable, then we can construct $W$ as an orientable manifold.
\end{Prop}

The following corollary has been shown.

\begin{Cor}[\cite{kitazawa 3}]
\label{cor:3}
In the situation of Proposition \ref{prop:8}, let $M$ be connected and $i:M \rightarrow W$ be
 the natural inclusion. Then \\
$$r_{\ast} \circ i_{\ast}:{\pi}_k(M) \rightarrow {\pi}_k(V)$$
 gives an isomorphism for $0 \leq k \leq m-\dim{V}-1$.
\end{Cor}

By using the results, the following has been shown.

\begin{Prop}[\cite{kitazawa 3}]
\label{prop:9}
Let $M$ be a closed and connected $C^{\infty}$ manifold of dimension $m$ and $f:M \rightarrow {\mathbb{R}}^n$ be a round fold
 map. Let $n \geq 2$ and $m-n \geq 2$. Suppose that $q_f(S(f))$ consists of two connected components and that one of the connected components is a $k$ S-locus. \\
\ \ \ We also assume that $f$ has a family of S-identifications. \\
\ \ \ Then there exist a compact ${\rm PL}$ manifold $W$ of dimension $m+1$ such that $\partial W=M$, a polyhedron $V$ and continuous maps $r:W \rightarrow V$
 and $s:V \rightarrow W_f$ and the followings hold.
\begin{enumerate}
\item $V$ is a polyhedron of dimension $n+k$. 
\item There exist a triangulation of $W$, a triangulation of $V$ and a triangulation of $W_f$ such that $r$ is a simplicial map and that $s$
 is a simplicial map and the followings hold.
\begin{enumerate}
\item For each $p \in V$, $r^{-1}(p)$ collapses to a point and $r$ is a homotopy equivalence.
\item If $p$ is in the closure of the AS-region in $W_f$, then $s^{-1}(p)$ is a point. If $p$ is in the AS-region in $W_f$ and
 $q \in s^{-1}(p)$, then
 $r^{-1}(q)$ is ${\rm PL}$ homeomorphic to $D^{m-n+1}$. 
\item If $p$ is in the $k$ S-region in $W_f$, then $s^{-1}(p)$ is ${\rm PL}$ homeomorphic to $S^k$. If $p$ is in the $k$ S-region
 in $W_f$ and $q \in s^{-1}(p)$, then $r^{-1}(q)$ is ${\rm PL}$ homeomorphic
 to $D^{m-n-k+1}$.
\end{enumerate}
\item $W$ collapses to a subpolyhedron ${V}^{\prime}$ such that $r {\mid}_{{V}^{\prime}}:{V}^{\prime} \rightarrow V$ is a ${\rm PL}$ homeomorphism.
\item $r_{\ast} \circ i_{\ast}:{\pi}_j(M) \rightarrow {\pi}_j(V)$ gives an isomorphism for $0 \leq j \leq m-\dim{V}-1=m-n-k-1$ where $i:M \rightarrow W$ is the inclusion.
\item $V$ is ${\rm PL}$ homeomorphic to $S^{n+k} {\bigcup}_{\psi} (S^{n-1} \times [0,1])$ for a ${\rm PL}$
 homeomorphism $\psi:B \rightarrow A$ {\rm (}a polyhedron obtained by glueing $S^{n+k}$ and $S^{n-1} \times [0,1]$ by $\psi${\rm )} where $A$ is ${\rm PL}$ homeomorphic to $S^{n-1}$ and
 trivially embedded in $S^{n+k}$ in the ${\rm PL}$ category and where $B:=S^{n-1} \times \{0\} \subset S^{n-1} \times [0,1]$.
\end{enumerate}
\end{Prop}

About Proposition \ref{prop:9} in the present paper, the following example also appeared in \cite{kitazawa 3}. 

\begin{Ex}[\cite{kitazawa 3}]
\label{ex:2}
It is also known that any closed and $2$-connected manifold of dimension $6$ has a $C^{\infty}$ differentiable structure and that the
 resulting $C^{\infty}$ manifold is always $C^{\infty}$ diffeomorphic to a connected sum of finite copies of $S^3 \times S^3$ (\cite{wall}). So in the situation
 of Proposition \ref{prop:9}, if $m=6$, $n=2$ and $k=1$, then $M$ is $C^{\infty}$ diffeomorphic to $S^3 \times S^3$. 
\end{Ex}

About the diffeomorphism types of the source manifolds of Proposition \ref{prop:9}, we show the following in the present paper.

\begin{Thm}
\label{thm:2}
For any pair of $C^{\infty}$ homotopy spheres ${\Sigma}_1$ and ${\Sigma}_2$ such that {\rm (}oriented{\rm )} ${\Sigma}_1$ is in a class of ${\Theta}_{m,n}$ {\rm (}$m>n \geq 2${\rm )} and that the
 dimension of ${\Sigma}_2$ is positive, any closed $C^{\infty}$ manifold $M$ having the structure of a $C^{\infty}$ ${\Sigma}_2$-bundle over ${\Sigma}_1$ admits a $C^{\infty}$ trivial round
 fold map into ${\mathbb{R}}^n$ satisfying the assumption of Proposition \ref{prop:9}.    
\end{Thm}

\begin{proof}
We apply the method used in the proof of Theorem \ref{thm:1} to a round fold map from ${\Sigma}_1$ into ${\mathbb{R}}^n$ whose singular set is connected. We obtain a $C^{\infty}$ trivial round fold map $f:M \rightarrow {\mathbb{R}}^n$ satisfying the followings.
\begin{enumerate}
\item The singular set $S(f)$ consists of two connected components.
\item The inverse image of an axis of the resulting map $f$ is $C^{\infty}$ diffeomorphic to $D^{m-n+1} \times {\Sigma}_2$.
\end{enumerate}
We easily know that $f$ satisfies the assumption of Proposition \ref{prop:9} by the proof of Theorem \ref{thm:1}.
\end{proof}

\section{P-operations and applications}
\label{sec:5}
\subsection{P-operations}
We easily have the following proposition.

\begin{Prop}
\label{prop:10}
Let $M$ be a closed $C^{\infty}$ manifold of dimension $m$ and $f:M \rightarrow {\mathbb{R}}^n$ be a locally $C^{\infty}$ trivial round
 fold map. Let $F$ be a closed $C^{\infty}$ manifold and ${M}^{\prime}$ be a closed $C^{\infty}$ manifold having the structure of a $C^{\infty}$
 $F$-manifold such that for any connected component $C$ of $f(S(f))$ and a small $C^{\infty}$ closed tubular neighborhood $N(C)$ of $C$ as in Definition \ref{def:8}, the
 restriction to $f^{-1}(N(C))$ is a trivial $C^{\infty}$ bundle. Then on ${M}^{\prime}$ there exists a locally $C^{\infty}$ trivial round
 fold map $f^{\prime}:{M}^{\prime} \rightarrow {\mathbb{R}}^n$. 
\end{Prop} 
\begin{proof}
In the proof of Proposition \ref{prop:7}, by replacing ${\phi}_k \times {\rm id}_F$ by a $C^{\infty}$
 diffeomorphism ${\Phi}_k$ which is regarded as a general isomorphism between the two trivial $C^{\infty}$ $F$-bundles inducing the
 identification ${\phi}_k$ of the base spaces, ${M}^{\prime}$, which has the structure of a $C^{\infty}$ $F$-bundle over $M$, is
 regarded as
 $$(\cdots ((f^{-1}({D^{n}}_{\frac{1}{2}}) \times F) {\bigcup}_{{\Phi}_1} (f^{-1}(P_1) \times F)) \cdots) {\bigcup}_{{\Phi}_l} ({f^{-1}}(P_l) \times F)$$
By analogy of the arguments of the proofs of Proposition \ref{prop:7} and Theorem \ref{thm:1}, we obtain a new round fold map $f^{\prime}:{M}^{\prime} \rightarrow {\mathbb{R}}^n$.
\end{proof}

In the proof, from $f:M \rightarrow {\mathbb{R}}^n$, we obtain ${f}^{\prime}:M^{\prime} \rightarrow {\mathbb{R}}^n$. We call the operation of
 constructing ${f}^{\prime}$ from $f$ a {\it P-operation by $F$ on $f$}. \\ 
\ \ \ In this paper, we consider $P$-operations by $S^1$. We introduce fundamental terms and known facts on $S^1$-bundles. It is well
 known that $C^{\infty}$ $S^1$-bundles are regarded as linear bundles whose structure groups are $O(2)$. A $C^{\infty}$ $S^1$-bundle is said
 to be {\it orientable} if the structure group is regarded as $SO(2)$. Since $C^{\infty}$ $S^1$-bundles are linear bundles, so we can consider
 their {\it Stiefel-Whitney classes}, which are cohomology classes of the base spaces with coefficients $\mathbb{Z}/2\mathbb{Z}$, including their {\it 1st Stiefel-Whitney classes} ({\it 2nd Stiefel-Whitney classes}), which are 1st cohomology classes (resp. 2nd cohomology classes) of the base spaces with coefficients $\mathbb{Z}/2\mathbb{Z}$. If a $C^{\infty}$ $S^1$-bundle is
 orientable, then we can orient the bundle and consider its {\it Euler class}, which is a 2nd cohomology class of the base space with
 coefficient $\mathbb{Z}$. We introduce known facts on classifications of linear $S^1$-bundles without proofs.

\begin{Prop}
\label{prop:11}
Let $X$ be a topological space.
\begin{enumerate}
\item The 1st Stiefel-Whitney class $\alpha \in H^1(X,\mathbb{Z}/2\mathbb{Z})$ of a linear $S^1$-bundle over $X$ vanishes if and only if the bundle is orientable.
\item For any $\alpha \in H^2(X,\mathbb{Z})$, there exists a linear oriented $S^1$-bundle over $X$ whose Euler class is $\alpha$. 
\item Two linear oriented $S^1$-bundles on $X$ are isomorphic if the Euler classes are same.
\end{enumerate}
\end{Prop}
 The facts are important in the following subsection. For general theory of linear (vector) bundles and their
 characteristic classes including Stiefel-Whitney classes and Euler classes, see also \cite{milnor stasheff} for example.

\subsection{Round fold maps on manifolds having the structures of $C^{\infty}$ $S^1$-bundles}
\begin{Thm}
\label{thm:3}
Let $M$ be a closed $C^{\infty}$ manifold of dimension $m$, $f:M \rightarrow {\mathbb{R}}^n$ be a 
locally $C^{\infty}$ trivial round fold map and $m \geq n \geq 3$. Assume that either of the followings holds.
\begin{enumerate}
\item $n \geq 4$.
\item $n=3$, $M$ is connected, $f(M)$ is $C^{\infty}$ diffeomorphic to $D^3$ and $W_f$ is connected.
\end{enumerate}
 Assume that for any connected component $C$ of $f(S(f))$, a normal fiber $F_C$ of $C$ corresponding to a trivial bundle over $C$ as in Definition \ref{def:8} satisfies $H^{2}(F_C,\mathbb{Z}) \cong \{0\}$. Then
 for any closed $C^{\infty}$
 manifold ${M}^{\prime}$ having the structure of a $C^{\infty}$ orientable $S^1$-bundle over $M$, by a P-operation by $S^1$ on $f$, we
 can obtain a round fold map $f^{\prime}:M^{\prime} \rightarrow {\mathbb{R}}^n$. \\ 
\ \ \ Furthermore, for any connected component $C$ of $f(S(f))$, let the fiber $F_C$ above satisfy $H^{1}(F_C,{\mathbb{Z}}_2) \cong \{0\}$. Then for any closed $C^{\infty}$
 manifold having the structure of a $C^{\infty}$ $S^1$-bundle over $M$, by a P-operation by $S^1$ on $f$, we
 can obtain a locally $C^{\infty}$ trivial round fold map $f^{\prime}:M^{\prime} \rightarrow {\mathbb{R}}^n$.
\end{Thm}
\begin{proof}
First we prove the case where $n \geq 4$. \\
\ \ \ For any connected component $C$ of $f(S(f))$, there exist a small $C^{\infty}$ closed tubular
 neighborhood $N(C)$ and a trivial $C^{\infty}$ bundle $f^{-1}(N(C))$ over $C$ as in Definition \ref{def:8} and $H^{2}(f^{-1}(N(C)),\mathbb{Z}) \cong \{0\}$ holds since $C$ is homeomorphic to $S^{n-1}$ and by the assumption, for a normal fiber $F_C$ of $C$ corresponding to $f^{-1}(N(C))$, we may assume that $H^{2}(F_C,\mathbb{Z}) \cong \{0\}$ holds. If we restrict
 any $C^{\infty}$ orientable $S^1$-bundle over $M$ to $f^{-1}(N(C))$, then the bundle obtained by the restriction is a trivial $C^{\infty}$ bundle. We
 can do a P-operation by $S^1$ on $f$. \\
\ \ \ Assume that $H^{1}(F_C,{\mathbb{Z}}/2\mathbb{Z}) \cong \{0\}$ holds. We note that $H^{1}(f^{-1}(N(C)),\mathbb{Z}/2\mathbb{Z}) \cong \{0\}$ holds since $C$ is homeomorphic to $S^{n-1}$
 and $H^{1}(F_C,\mathbb{Z}/2\mathbb{Z}) \cong \{0\}$ holds. If we restrict any $C^{\infty}$ $S^1$-bundle over $M$ to $f^{-1}(N(C))$, then the bundle obtained by the restriction is orientable and as a result it is a trivial $C^{\infty}$ bundle. We can do a P-operation similarly. \\
\ \ \ Thus this completes the proof of the case where $n \geq 4$. Second we prove the case where $n=3$, where $f(M)$ is $C^{\infty}$ diffeomorphic to $D^3$ and where $W_f$ is connected. For any connected component $C$ of $f(S(f))$, there
 exist a small $C^{\infty}$ closed tubular
 neighborhood $N(C)$ and a trivial $C^{\infty}$ bundle $f^{-1}(N(C))$ over $C$ as in Definition \ref{def:8} and $H^{2}(f^{-1}(N(C)),\mathbb{Z}) \cong H^{2}(C,\mathbb{Z}) \cong \mathbb{Z}$ since $C$
 is homeomorphic to $S^2$ and by the assumption, for a normal fiber $F_C$ of $C$ corresponding to the bundle $f^{-1}(N(C))$, we may assume that $H^{2}(F_C,\mathbb{Z}) \cong \{0\}$ holds. \\
\ \ \ From the cohomology group of the fiber $F_C$, it suffices to show 
 that the restriction of any $C^{\infty}$ orientable $S^1$-bundle over $M$ to $f^{-1}(\partial N(C))$ is a trivial $C^{\infty}$ bundle for any $C$. For this, it suffices to
 show that the Euler class of the restriction of the (oriented) bundle to $f^{-1}(\partial N(C))$ vanishes. \\
\ \ \ For $C$ and a connected component $C_1$ of $\partial N(C)$, if the restriction of
 the bundle over $M$ to $f^{-1}(C_1)$ is orientable and the Euler class of the obtained (oriented) bundle vanishes, then for the other connected component $C_2$ of
 $\partial N(C)$, the restriction of the bundle over $M$ to $f^{-1}(C_2)$ is orientable and the Euler class of the obtained (oriented) bundle vanishes. \\
\ \ \ Let $F$ be the fiber of a point in a proper core $P$ of $f$. Note that for a proper core $P$ of $f$, $f^{-1}(P)$ has the structure
 of a trivial $C^{\infty}$ $F$-bundle over $P$. $f^{-1}(\partial P)$ is regarded as the restriction of
 the bundle $f^{-1}(P)$ to $\partial P$ and for the connected component $C_0$ of $f(S(f))$ in the center
 and the given small $C^{\infty}$ closed tubular neighborhood $N(C_0)$, it is regarded as a subbundle of the trivial $C^{\infty}$ bundle $f^{-1}(N(C_0))$ over $C_0$, whose fiber
 is $F_{C_0}$ satisfying $H^{2}(F_{C_0},\mathbb{Z}) \cong \{0\}$. It means that the Euler
 class of the restriction of the bundle over $M$ to any connected
 component of $f^{-1}(\partial P)$ vanishes if oriented. \\
\ \ \ Since $W_f$ is connected, for any connected comonent $C$ of $f(S(f))$, it follows that the Euler class of the restriction of a $C^{\infty}$ orientable $S^1$-bundle over $M$ to $f^{-1}(\partial N(C))$ vanishes if oriented. \\
\ \ \ Assume that $H^{1}(F_C,{\mathbb{Z}}/2\mathbb{Z}) \cong \{0\}$ also holds. We note that $H^{1}(f^{-1}(N(C)),\mathbb{Z}/2\mathbb{Z}) \cong \{0\}$ holds since $C$ is homeomorphic to $S^{n-1}$
 and $H^{1}(F_C,\mathbb{Z}/2\mathbb{Z}) \cong \{0\}$ holds. If we restrict any $C^{\infty}$ $S^1$-bundle over $M$ to $f^{-1}(N(C))$, then the bundle obtained by the restriction is orientable
 and as a result it is a trivial $C^{\infty}$ $S^1$-bundle. \\
\ \ \ Thus this completes the proof of the theorem.
\end{proof}

\begin{figure}
\begin{center}
\includegraphics[width=50mm]{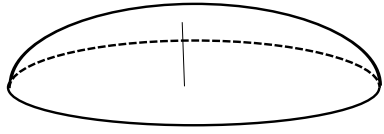}
\end{center}
\caption{A fiber of the bundle $V_C$. (The segment in this disc is a fiber of the bundle $\widetilde{P(C)}$ and the submanifold of the boundary
 above the corner is a fiber of the bundle ${q_f}^{-1}(N(C))$.)}
\label{fig:3}
\end{figure}

  For a closed $C^{\infty}$ manifold $M$ of dimension $m$, let $f:M \rightarrow {\mathbb{R}}^n$ ($m \geq n \geq 2$) be
 a round fold map which is locally $C^{\infty}$ trivial and which has a
 family of S-identifications. We also assume that $n \geq 3$ and that for any integer $k$ satisfying $m-n-k<3$, there exist no $k$-loci in $W_f$. Last we also
 assume that $W_f$ has no branched points. \\
\ \ \ First, for $f$, we demonstrate the proof of Proposition \ref{prop:8} of the present paper or Theorem \ref{thm:1} of \cite{kitazawa 3}. \\
\ \ \ Let $C$ be a connected component of the image of the set $F_0(f)$ of all the definite fold points of $f$. Let $N(C)$ be a small
 regular neighborhood of $C$. We note that $N(C)$ has the structure of a trivial ${\rm PL}$ bundle over $C$ with a fiber homeomorphic to $[0,1]$. We may assume that $C$ corresponds to the $0$-section ($0 \in [0,1]$). \\
\ \ \ For each $p \in C$, set $K_p:={q_f}^{-1}(\{p\} \times [0,1])$ for a fiber $\{p\} \times [0,1]$ of the bundle $N(C)$ over $C$. It
 is $C^{\infty}$ diffeomorphic to $D^{m-n+1}$. We may assume that $q_f^{-1}(N(C))$ has the structure
 of a $C^{\infty}$ bundle over $C$ with a fiber $C^{\infty}$ diffeomorphic to $D^{m-n+1}$ and $K_p$ is the fiber of $p \in C$. 
Now we can construct a $C^{\infty}$ compact ($m+1$)-manifold $V_C$ having the structure of a $C^{\infty}$ $D^{m-n+2}$-bundle
 over $C$ (we denote by $\widetilde{K_p}$ the fiber of $p \in C$) such that $q_f^{-1}(N(C))$ has the structure of
 a subbundle of the bundle $V_C$ with a fiber $C^{\infty}$ diffeomorphic to $D^{m-n+1}$ and is in the
 boundary of $V_C$ and a ${\rm PL}$ map $r_C:V_C \rightarrow N(C)$
such that ${r_C}^{-1}(p,t)$ is ${\rm PL}$ homeomorphic to $D^{m-n+1}$ for $p \in C$ and $t \in (0,1]$ and
 that ${r_C}^{-1}(p,t)$ is a point for $p \in C$ and $t=0$. Furthermore, there exists
 a ${\rm PL}$ manifold $\widetilde{P(C)} \subset V_C$ of dimension $n$ such that $\widetilde{P(C)}$
 has the structure of a subbundle of the bundle $V_C$, that ${r_C} {\mid}_{\widetilde{P(C)}}:\widetilde{P(C)} \rightarrow N(C)$ is
 a ${\rm PL}$ homeomorphism (a bundle isomorphism between the two ${\rm PL}$ bundles), that
$\widetilde{P(C)} \bigcap \partial \widetilde{K_p}$ consists of a point $(p,1) \in \{p\} \times [0,1]$ and is not
 in $q_f^{-1}(N(C))$ and that $V_C$ collapses to $\widetilde{P(C)}$. See also FIGURE \ref{fig:3}. 
Finally we set $s_C:={\rm id}_{N(C)}$ for latter discussions. \\
\ \ \ Let a connected component $C$ of $q_f(S(f))$ be a $k-1$ S-locus where $2 \leq k<m-n$. \\
\ \ \ By the assumption, we can take a regular neighborhood $N(C)$ so that ${q_f}^{-1}(N(C))$ has the structure
 of a trivial $C^{\infty}$ bundle over $C$ with a fiber ${\rm PL}$ homeomorphic to $D^k \times S^{m-n-k+1}-{\rm Int} D^{m-n+1}$.
 We note that $N(C)$
 has the structure of a trivial
 ${\rm PL}$ $[-1,1]$-bundle over $C$ and we may assume that $C$ corresponds to the $0$-section ($0 \in [-1,1]$). We may assume that there exists
 a ${\rm PL}$ bundle over $C$ with a fiber ${\rm PL}$ homeomorphic to $(S^{m-n-k+1} \times D^k) {\bigcup}_{{\phi}_1 \times {\phi}_2} (D^{m-n-k+2} \times S^{k-1})$ (and $S^{m-n+1}$) for a pair
 of ${\rm PL}$ homeomorphisms
 $({\phi}_1:\partial D^{m-n-k+2} \rightarrow S^{m-n-k+1}, {\phi}_2:S^{k-1} \rightarrow \partial D^k)$ and that the
 bundle ${q_f}^{-1}(N(C))$ is a subbundle of the
 $(S^{m-n-k+1} \times D^k)$-bundle. For this fact, remember
 that by the definition of S-identifications, for any connected component ${C}^{\prime}$ of $\partial N(C)$ in the $k-1$ S-region, ${q_f}^{-1}({C}^{\prime})$
 has the structure of a subbundle of the bundle ${q_f}^{-1}(N(C))$ over ${C}^{\prime}$ with a fiber ${\rm PL}$ homeomorphic
 to $\partial D^k \times S^{m-n-k+1}$ and that the structure group of ${q_f}^{-1}({C}^{\prime})$ is a subgroup of
 the $C^{\infty}$ diffeomorphism group of the fiber, consisting of some $C^{\infty}$ diffeomorphisms regarded as bundle isomorphisms on the
 trivial ${\rm PL}$ bundle $\partial D^k \times S^{m-n-k+1}$ over $\partial D^k$ inducing ${\rm PL}$ homeomorphisms on
 the base space $\partial D^k=S^{k-1}$. \\
\ \ \ Then we can construct the following objects.

\begin{enumerate}
\item A compact ${\rm PL}$ manifold $V_C$ of dimension $m+1$ having the structure of
 a trivial ${\rm PL}$ $D^{m-n+2}$-bundle over $C$ (let ${V_C}_p$ denote the
 fiber of $p$ of the base space) such that the subbundle corresponding to
 $\partial D^{m-n+2}=(S^{m-n-k+1} \times D^k) {\bigcup}_{{\phi}_1 \times {\phi}_2} (D^{m-n-k+2} \times S^{k-1})$ is the ${\rm PL}$ $S^{m-n+1}$-bundle in the previous paragraph. The bundle $q_f^{-1}(N(C))$
 (over $C$) is
 a subbundle of the bundle $V_C$ and is in the boundary of $V_C$ (we regard it is also a subbundle of the ($S^{m-n-k+1} \times D^k$)-bundle).
\item A trivial ${\rm PL}$ bundle $P(C)$ whose fiber is the mapping cylinder of the constant
 map $c:S^{k-1} \rightarrow [0,1]$ satisfying $c(S^{k-1})=\{0\}$ and whose base space
 is $C$.
\item ${\rm PL}$ maps $r_C:V_C \rightarrow P(C)$ and $s_C:P(C) \rightarrow N(C)$
 such that ${s_C}^{-1}(p)$ is a point for $p$ in the closure of an AS-region, that ${s_C}^{-1}(p)$ is ${\rm PL}$ homeomorphic to $S^{k-1}$
 for $p$ in a $k-1$ S-region, that ${r_C}^{-1}(q)$ is ${\rm PL}$ homeomorphic to $D^{m-n+1}$ for $q \in {s_C}^{-1}(p)$ for $p$ in an AS-region,
 that ${r_C}^{-1}(q)$ is ${\rm PL}$ homeomorphic to $D^{m-n-k+2}$ for $q \in {s_C}^{-1}(p)$ for $p$ in a $k-1$ S-region and that ${(s_C \circ r_C)}^{-1}(\{p\} \times [-1,1])={V_C}_p$ for all $p \in C$ ($\{p\} \times [-1,1]$ is
 a fiber of the bundle $N(C)$ over $C$; $N(C)$ has the structure of a
 trivial ${\rm PL}$ $[-1,1]$-bundle over $C$ and $C$ corresponds to the $0$-section). Furthermore, ${r_C}^{-1}(q)$
 collapses to a point for $q \in P(C)$. 
\item A subpolyhedron $\widetilde{P(C)} \subset V_C$ of dimension $n+k-1$ such
 that $\widetilde{P(C)}$ has the structure of a subbundle of the trivial bundle $V_C$,
 that ${r_C} {\mid}_{\widetilde{P(C)}}:\widetilde{P(C)} \rightarrow P(C)$ is
 a ${\rm PL}$ homeomorphism (a bundle isomorphism between the two ${\rm PL}$ bundles), that $\widetilde{P(C)} \bigcap \partial {V_C}_p=\{p\} \times ((S^{k-1} \times \{0\}) \sqcup \{1\}) \subset \widetilde{P(C)}$, that $\widetilde{P(C)} \bigcap \partial {V_C}_p$ does not contain any points in ${q_f}^{-1}(N(C))$ and is the disjoint union of a point in $S^{m-n-k+1} \times D^{k}$ and a ($k-1$)-sphere $\{0\} \times S^{k-1}$,
 and that $V_C$ collapses to $\widetilde{P(C)}$.
\end{enumerate}     

\ \ \ In this argument, we may assume that $f$ is a normal form of a round fold map and that for any connected component $C$ of $f(S(f))$ and
 a positive integer $j$, $N(C)=[j-\frac{1}{2},j+\frac{1}{2}]$. \\
\ \ \ Let $f(M)$ be $C^{\infty}$ diffeomorphic to $D^n$. Then, $f {\mid}_{f^{-1}({D^n}_{\frac{1}{2}})} \rightarrow {D^n}_{\frac{1}{2}}$ gives the structure of
 a trivial $C^{\infty}$ bundle. $D:={\bar{f}}^{-1}({D^n}_{\frac{1}{2}}) \subset W_f$ is in a connected component of $W_f-q_f(S(f))$. \\
\ \ \ Let the fiber of a point in $D$ be an almost-sphere and ${q_f}^{-1}(D)$ have the structure of
 a $C^{\infty}$ bundle over $D$ with a fiber ${\rm PL}$ homeomorphic to $S^{m-n}$. Let
 $r_D:V_D \rightarrow D$ give the structure of a ${\rm PL}$ $D^{m-n+1}$-bundle which is an associated bundle of
 the bundle ${q_f}^{-1}(D)$ over $D$ (we regard $S^{m-n}=\partial D^{m-n+1}$). We take the associated bundle so that the structure group
 is a group consisting of
 ${\rm PL}$ homeomorphisms $r$ on $D^{m-n+1}$ such that $r(0)=0$ and that for a ${\rm PL}$
 homeomorphism ${r}^{\prime}$ on $S^{m-n}$, $\frac{r(x)}{|x|}={r}^{\prime}(\frac{x}{|x|})$ ($x \neq 0$). Let $\widetilde{P(D)} \subset V_D$ be the
 $0$-section of the associated bundle. Finally we set $P(D):=D$ and let $s_{D}:P(D) \rightarrow D$ be the identity map. \\
\ \ \ Let $D$ be in a $k$ S-region whose closure is bounded
 by a $k$ S-locus and ${q_f}^{-1}(D)$ have the structure of
 a $C^{\infty}$ bundle over $D$ with a fiber ${\rm PL}$ homeomorphic to $S^k \times S^{m-n-k}$. Let $V_{D}$ have the structure of a ${\rm PL}$ ($S^k \times D^{m-n-k+1}$)-bundle over $D$ which is
 an associated bundle of the bundle ${q_f}^{-1}(D)$ over $D$ (we regard $S^k \times S^{m-n-k}=S^k \times \partial D^{m-n-k+1}$). We take the associated bundle so that the structure group
 is a group consisting of
 some ${\rm PL}$ homeomorphisms $r$ on $S^k \times D^{m-n-k+1}$ such that $r$ is a bundle isomorphism on
 a ${\rm PL}$ bundle $S^k \times D^{m-n-k+1}$ over $S^k$ whose structure group consists of ${\rm PL}$ homeomorphisms ${r_1}$ on $D^{m-n-k+1}$ where $r_1(0)=0$ and where for a ${\rm PL}$ 
homeomorphism ${r_2}$ on $S^{m-n-k}$, $\frac{r_1(x)}{|x|}={r_2}(\frac{x}{|x|})$ ($x \neq 0$) and that $r$ induces a ${\rm PL}$ homeomorphism on the base space $S^k$. Let $\widetilde{P(D)} \subset V_D$ be the
 subbundle whose fiber is $S^k \times \{0\} \subset S^k \times D^{m-n-k+1}$. The bundle $V_D$ over $D$ is also given by the
 composition of two ${\rm PL}$ maps $r_{D}:V_D \rightarrow P(D)$ and $s_D:P(D) \rightarrow D$ where
 $P(D)$ is a ${\rm PL}$ manifold and where $s_D$ gives $P(D)$ the structure of a bundle over $D$ equivalent to the bundle $\widetilde{P(D)}$ over $D$. \\
\ \ \ Now, by the definition
 of a family of S-identifications, we obtain $W$, $V$, $V^{\prime}$, $r$ and $s$ in
 Proposition \ref{prop:8} by glueing a family $\{V_C \}$ and $V_D$, a family $\{P(C) \}$ and $P(D)$, a family $\{\widetilde{P(C)} \}$ and $\widetilde{P(D)}$, a family $\{r_C \}$ and $r_D$ and a family $\{s_C \}$ and $s_D$ together, respectively. This completes
 the proof of Proposition \ref{prop:8} for $f$. If $f(M)$ is not $C^{\infty}$ diffeomorphic to $D^n$, then we don't need $V_D$, $P(D)$, $\widetilde{P(D)}$, $r_D$ and $s_D$. \\   
\ \ \ We also note that for any axis $L$ of $f$ and any $k$ S-region $R$ in $W_f$ whose closure $\overline{R}$ is bounded by a disjoint
 union of a pair of two $k$ S-loci, ${s \circ \bar{f}}^{-1}(\overline{R} \bigcap L)$ is ${\rm PL}$ homeomorphic to $S^{k+1}$. Thus the simple
 homotopy type of $V$ is represented as in the following paragraph. \\
\ \ \ Let $A:=D^n$ if the fiber of a point in a proper core
 of $f$ is an almost-sphere and $A:=S^{n+k}$ if $f(M)$ is $C^{\infty}$ diffeomorphic to $D^n$ and the connected component
 of $W_f-q_f(S(f))$ homeomorphic to an open disc of dimension $n$ is a $k$ S-region whose closure is bounded by a $k$ S-locus. Let $l_k$ be the number of S-regions bounded by disjoint unions of two $k$ S-loci in $W_f$. Let $L$ be a
 wedge product of $l_k$ copies of $S^{k+1}$ for all $1 \leq k \leq m-n-1$. If $f(M)$ is $C^{\infty}$ diffeomorphic to $D^n$, then we define
 a ${\rm PL}$ homeomorphism $\phi$ from $S^{n-1} \times \{\ast\} \subset B:=S^{n-1} \times L$ onto $\partial A$ in the case where $A=D^n$ and onto a standard sphere of
 dimenison $n-1$ trivially embedded in $A$ in the ${\rm PL}$ category in the case where $A=S^{n+k}$. Hence $V$ in the proposition is simple homotopy
 equivalent to $A {\bigcup}_{\phi} B$ if $f(M)$ is $C^{\infty}$ diffeomorphic to $D^n$ and $V$ is simple homotopy equivalent to $B=S^{n-1} \times L$ if $f(M)$ is not $C^{\infty}$ diffeomorphic to $D^n$. \\ 
\ \ \ Furthermore, ${\pi}_1(V) \cong {\pi}_1(M) \cong H_1(M,\mathbb{Z}) \cong H^1(M,\mathbb{Z}) \cong H^1(M,\mathbb{Z}/2\mathbb{Z}) \cong \{0\}$ holds by Proposition \ref{prop:8} and Corollary \ref{cor:3}. By these propositions, ${\pi}_2(V) \cong {\pi}_2(M) \cong H_2(M,\mathbb{Z}) \cong H^2(M,\mathbb{Z}) \cong {\mathbb{Z}}^{l_1}$ also holds
 where $l_1$ is the number of $1$ S-regions in $W_f$ whose closures are bounded by pairs of two $k$ S-loci defined above. \\
\ \ \ $f$ satisfies the assumption of Theorem \ref{thm:3}. Now this completes the proof of the following theorem.

\begin{Thm}
\label{thm:4}
Let $M$ be a closed $C^{\infty}$ manifold of dimension $m$. Let $f:M \rightarrow {\mathbb{R}}^n$ be a round
 fold map which is locally $C^{\infty}$ trivial and which has a
 family of S-identifications. Assume that the Reeb space $W_f$ has no branched points. Assume also
 that $n \geq 3$ and that for any integer $k$ satisfying $m-n-k<3$, there exist no $k$-loci in $W_f$. Suppose that either of the followings holds.
\begin{enumerate}
\item $n \geq 4$.
\item $n=3$, $M$ is connected, $f(M)$ is $C^{\infty}$ diffeomorphic to $D^3$ and $W_f$ is connected.
\end{enumerate}
Let $l_k$ be the number of S-regions bounded by disjoint unions of pairs of two $k$ S-loci in $W_f$. Then $H^2(M,\mathbb{Z})$ is a
 free commutative group of rank $l_1$ and on a manifold having the structure of a $C^{\infty}$ oriented $S^1$-bundle whose
 Euler class is $\alpha \in H^2(M,\mathbb{Z})$, we can obtain a round fold map by a P-operation by $S^1$ on $f$.  
\end{Thm}

We have the following theorem where the target manifold is ${\mathbb{R}}^2$.

\begin{Thm}
\label{thm:5}
Let $M$ be a closed $C^{\infty}$ manifold of dimension $m (m \geq 2)$ and $f:M \rightarrow {\mathbb{R}}^2$ be a 
locally $C^{\infty}$ trivial round fold map.
\begin{enumerate}
\item For any connected component $C$ of $f(S(f))$, we denote a small $C^{\infty}$ closed tubular
 neighborhood of $C$ as in Definition \ref{def:8} by $N(C)$ and we denote a normal fiber of $C$ corresponding to the trivial
 bundle $f^{-1}(N(C))$ over $C$ as in Definiiton \ref{def:8} by $F_C$. Let $F_C$ satisfy $H^{1}(F_C,\mathbb{Z}) \cong H^{2}(F_C,\mathbb{Z}) \cong \{0\}$. Then by a P-operation by $S^1$ on $f$, for any closed $C^{\infty}$
 manifold having the structure of a $C^{\infty}$ $S^1$-bundle over $M$ such that
 for any connected component $C$ of $f(S(f))$, the restriction to $f^{-1}(\partial N(C))$ is an orientable bundle, we
 can obtain a round fold map $f^{\prime}:M^{\prime} \rightarrow {\mathbb{R}}^2$. \\
\item Furthermore, if $f(M)$ is $C^{\infty}$ diffeomorphic
 to $D^2$ and $W_f$ is connected, then by a P-operation by $S^1$ on $f$, for any closed $C^{\infty}$
 manifold $M^{\prime}$ having the structure of a $C^{\infty}$ $S^1$-bundle over $M$, we
 can obtain a round fold map $f^{\prime}:M^{\prime} \rightarrow {\mathbb{R}}^2$.
\end{enumerate}
\end{Thm}
\begin{proof}
We prove the first part. For any connected component $C$ of $f(S(f))$, $H^{2}(f^{-1}(N(C)),\mathbb{Z}) \cong \{0\}$ holds since $C$ is homeomorphic to $S^1$ and
 for a normal fiber $F_C$ of $C$ corresponding to the bundle $f^{-1}(N(C))$ over $C$, $H^{1}(F_C,\mathbb{Z}) \cong H^{2}(F_C,\mathbb{Z}) \cong \{0\}$ holds. Note
 that $H^{1}(F_C,\mathbb{Z}/2\mathbb{Z}) \cong \{0\}$ holds. If we
 restrict any $C^{\infty}$ $S^1$-bundle over $M$ to $f^{-1}(N(C))$ such that
 the restriction to $f^{-1}(\partial N(C))$ is a $C^{\infty}$ orientable $S^1$-bundle, then the bundle obtained by the restriction is orientable since its 1st Stiefel-Whitney class vanishes and as a result, it is a
 trivial $C^{\infty}$ bundle. By a P-operation by $S^1$ on $f$, a round fold map $f^{\prime}:M^{\prime} \rightarrow {\mathbb{R}}^n$ is
 obtained. \\ 
\ \ \ This completes the proof of the first part. Now we prove the second part. For any connected component $C$ of $f(S(f))$, $H^{1}(f^{-1}(N(C)),\mathbb{Z}/2\mathbb{Z}) \cong H^{1}(C,\mathbb{Z}/2\mathbb{Z}) \cong \mathbb{Z}/2\mathbb{Z}$
 and $H^{2}(f^{-1}(N(C)),\mathbb{Z}) \cong {0}$ hold since $C$
 is homeomorphic to $S^1$ and $H^{1}(F_C,\mathbb{Z}) \cong H^1(F_C,\mathbb{Z}/2\mathbb{Z}) \cong H^{2}(F_C,\mathbb{Z}) \cong \{0\}$ holds. \\
\ \ \ From the cohomology group of the fiber $F_C$, it suffices to show 
 that the restriction of every $C^{\infty}$ $S^1$-bundle over $M$ to $f^{-1}(\partial N(C))$ is orientable and as a result, it is a trivial $C^{\infty}$ bundle. For this, it suffices to
 show that the 1st Stiefel-Whitney class of the restriction of the bundle over $M$ to $f^{-1}(\partial N(C))$ vanishes and that the Euler
 class of the restriction, which is defined by orienting the obtained bundle, vanishes. \\
\ \ \ For any connected component $C$ of $f(S(f))$ and a connected component $C_1$ of $\partial N(C)$, if the 1st
 Stiefel-Whitney class of the restriction of
 the $S^1$-bundle over $M$ to $f^{-1}(C_1)$ vanishes and the Euler
 class of the restrcition, which is defined by orienting the obtained bundle, vanishes, then for the other connected component $C_2$ of
 $\partial N(C)$, the 1st Stiefel-Whitney class of the restriction of
 the bundle over $M$ to $f^{-1}(C_2)$ vanishes and the Euler
 class of the restriction, which is defined by orienting the obtained bundle, vanishes. \\
\ \ \ Let $F$ be the fiber of a point in a proper core $P$ of $f$. Note that for a proper core $P$ of $f$, $f^{-1}(P)$ has the structure
 of a trivial $C^{\infty}$ $F$-bundle over $P$. $f^{-1}(\partial P)$ is regarded as the restriction of
 the bundle $f^{-1}(P)$ to $\partial P$ and for the connected component $C_0$ of $f(S(f))$ in the center
 and the given small $C^{\infty}$ closed tubular neighborhood $N(C_0)$, it is regarded as a subbundle of the trivial $C^{\infty}$ bundle $f^{-1}(N(C_0))$ whose fiber
 is $F_{C_0}$ satisfying $H^{1}(F_{C_0},\mathbb{Z}) \cong H^{1}(F_{C_0},\mathbb{Z}/2\mathbb{Z}) \cong H^{2}(F_{C_0},\mathbb{Z}) \cong \{0\}$. It means
 that the 1st Stiefel-Whitney class of the restriction of any $C^{\infty}$ $S^1$-bundle over $M$ to $f^{-1}(\partial P)$ vanishes and that the Euler
 class of the restriction of the bundle over $M$, which is defined by orienting the obtained bundle, vanishes. \\
\ \ \ Since $W_f$ is connected, by the arguments above, for any $C$, the restriction of any $C^{\infty}$ $S^1$-bundle over $M$ to $f^{-1}(\partial N(C))$ is orientable and the Euler class of the resulting (oriented) bundle vanishes. \\ 
\ \ \ Thus this completes the second part of the proof. This completes the proof. 
%
\end{proof}

\begin{Ex}
\label{ex:3}
A round fold map $f$ in Proposition \ref{prop:6} satisfies the assumptions of two statements of Theorem \ref{thm:5} for the case where $n=2$ and $m-n=m-2>2$. 
\end{Ex}

We review Theorem 3 of \cite{kitazawa 2} or Proposition 5 of \cite{kitazawa 4} (see also \cite{kitazawa}). 

\begin{Prop}[\cite{kitazawa}, \cite{kitazawa 2}]
\label{prop:12}
 Let $M$ be a closed and connected $C^{\infty}$ manifold of dimension $m$, $f:M \rightarrow {\mathbb{R}}^n$ be
 a round fold map and $m>n \geq 2$. If $m-n=1$, then we also assume that $M$ is orientable. \\
\ \ \ We assume that $f^{-1}(p)$ is a disjoint union of almost-spheres for each regular value $p$ and that the indices
 of all the fold points of $f$ are $0$ or $1$. \\     
\ \ \ Let $L$ be an axis of $f$ and $f_L:=f {\mid}_{f^{-1}(L)}$. We denote by $l_1$ the number of loops of the Reeb space $W_{f_{L}}$ of $f_L$ {\rm (}in
 other words, let $H_1(W_{f_L}) \cong {\mathbb{Z}}^{l_1}${\rm )}. We
 denote by $l_2$ the number of connected components of the fiber of a point in a proper core of $f$. \\  
\ \ \ Then there exist $W$ such that $\partial W=M$ and a homotopy equivalence $r:W \rightarrow W_f$. Furthermore, for the inclusion $i:M \rightarrow W$, $q_f=r \circ i$ gives
 an isomorphism of homotopy groups ${\pi}_{k}(M) \cong {\pi}_{k}(W_f)$ for $0 \leq k \leq m-n-1$ and we have the following list where
 we denote the free group of rank $r$ by $F_r$.

\begin{enumerate} 
\item  When $n \geq 3$, $m \geq 2n$, $f$ is topologically quasi-trivial and we have the followings.
$${\pi}_k(M) \cong {\pi}_k(W_f) \cong \\
\begin{cases}
F_{l_1} & k=1 \\
\{0\} & 2 \leq k < n-1
\end{cases}
$$

$${\pi}_{n-1}(M) \cong {\pi}_{n-1}(W_f) \cong \\
\begin{cases}
\mathbb{Z} & l_2=0 \\  
\{0\} & l_2 \neq 0
\end{cases}
$$
\item When $n \geq 3$, $n<m \leq 2n-1$ and $m-n \geq 2$, $f$ is topologically quasi-trivial and we have the following. 
$${\pi}_k(M) \cong {\pi}_k(W_f) \cong \\
\begin{cases}
F_{l_1} & k=1 \\
\{0\} & 2 \leq k \leq m-n-1
\end{cases}
$$
\item When $m \geq 4$ and $n=2$, we have the followings.
\begin{enumerate}
\item If $f$ is topologically quasi-trivial and $l_2=0$, then we have the following.
$${\pi}_k(M) \cong {\pi}_k(W_f) \cong \\
\begin{cases}
\mathbb{Z} \times F_{l_1} & k=1 \\
\{0\} & 2 \leq k \leq m-3
\end{cases}
$$
\item If $f$ is topologically quasi-trivial and $l_2 \neq 0$, then we have the following.
$${\pi}_1(M) \cong {\pi}_1(W_f) \cong F_{l_1}$$
\end{enumerate}
\end{enumerate}
\end{Prop}

If in Proposition \ref{prop:6} or Example \ref{ex:3}, $n=2$ and $m-n=m-2 \geq 3$, then ${\pi}_1(M) \cong {\pi}_1(W_f) \cong \{0\}$ and $H_2(M,\mathbb{Z}) \cong {\pi}_2(M) \cong {\pi}_2(W_f) \cong {\mathbb{Z}}^l$ hold by Proposition \ref{prop:12}
 where $l$ is the number appeared in Proposition \ref{prop:6}. As
 a result, $H^2(M,\mathbb{Z}) \cong {\mathbb{Z}}^{l}$ holds. For any class $\alpha \in H^2(M,\mathbb{Z})$, there exists a closed manifold $M(\alpha)$ having the structure of a $C^{\infty}$ oriented $S^1$-bundle over $M$ whose Euler class is $\alpha$ and on $M(\alpha)$ we have
 a round fold map into ${\mathbb{R}}^2$ by a P-operation by $S^1$ on the original map.

\begin{Ex}
\label{ex:4}
Let $M=S^2 \times S^3$. Any closed manifold ${M}^{\prime}$ having the structure of a $C^{\infty}$ orientable $S^1$-bundle over $M$ is represented as a product of a closed manifold having
 the structure of a $C^{\infty}$ orientable $S^1$-bundle over $S^2$ and $S^3$. We have a new round fold map into ${\mathbb{R}}^2$ on ${M}^{\prime}$ by a P-operation by $S^1$ on a round fold map $f:M \rightarrow {\mathbb{R}}^2$ as in Proposition \ref{prop:6} and Example \ref{ex:3}. \\ 
\ \ \ The family of all the manifolds having the structures of $C^{\infty}$ orientable $S^1$-bundles over $M$ contains a manifold $C^{\infty}$ diffeomorphic to $S^3 \times S^3$, for
 example. This gives a round fold map from $S^3 \times S^3$ into ${\mathbb{R}}^2$ and the map is different from one in Example 6 of \cite{kitazawa 3} or Example \ref{ex:2} of the present paper.   
\end{Ex}

\end{document}